\documentclass[final]{siamltex}

\usepackage{amsmath} 
\usepackage{amsfonts}

\usepackage{float}
\usepackage{graphicx}
\usepackage{subcaption}
\usepackage{url}
\usepackage{color}
\usepackage{bm}
\usepackage{multirow}
\usepackage{bigints}

\usepackage{listings}
\usepackage{color}

\usepackage{tikz}
\usetikzlibrary{trees, cd, babel}

\usepackage{todonotes}
\presetkeys%
    {todonotes}%
    {inline}{}
\definecolor{seabornblue}{rgb}{0.2980392156862745, 0.4470588235294118, 0.6901960784313725}  

\usepackage{hyperref}
\hypersetup{
    colorlinks,
    linkcolor={red!50!black},
    citecolor={blue!50!black},
    urlcolor={blue!80!black}
}

\title{An augmented Lagrangian preconditioner for the control of the Navier--Stokes equations}

\author{Santolo Leveque\footnotemark[1]
        \and Michele Benzi\footnotemark[2]
        \and Patrick E.~Farrell\footnotemark[3]}

\begin{document}
\maketitle
\renewcommand{\thefootnote}{\fnsymbol{footnote}}

\footnotetext[1]{CRM Ennio De Giorgi, Scuola Normale Superiore, Piazza dei Cavalieri 3, 56126, Pisa, Italy ({\tt santolo.leveque@sns.it}).}
\footnotetext[2]{Scuola Normale Superiore, Piazza dei Cavalieri 7, 56126 Pisa, Italy ({\tt michele.benzi@sns.it}).}
\footnotetext[3]{Mathematical Institute, University of Oxford, Radcliffe Observatory Quarter, Woodstock Road, Oxford OX2 6GG, 
UK ({\tt patrick.farrell@maths.ox.ac.uk}).}

\renewcommand{\thefootnote}{\arabic{footnote}}

\begin{abstract}
We address the solution of the distributed control problem for the steady,
	incompressible Navier--Stokes equations. We propose an inexact Newton
	linearization of the optimality conditions.
	Upon discretization by a finite element scheme, we obtain a sequence of
	large symmetric linear systems of saddle-point type. We use an
	augmented Lagrangian-based block triangular preconditioner in combination
	with the flexible GMRES method at each Newton step. The
	preconditioner is applied inexactly via a suitable multigrid solver.
	Numerical experiments indicate that the resulting method appears to be
	fairly robust with respect to viscosity, mesh size, and the choice of
	regularization parameter when applied to 2D problems.
\end{abstract}



\pagestyle{myheadings}
\thispagestyle{plain}
\markboth{}{AL PRECONDITIONER FOR NAVIER--STOKES CONTROL}

\section{Introduction}\label{sec_1}
In this work, we consider the distributed control of the flow of a Newtonian
	viscous fluid, subject to the incompressibility constraint. We restrict
	ourselves to steady problems; in this case, the physics
	is described by the stationary, incompressible Navier--Stokes equations.
	
The distributed control of the stationary incompressible Navier--Stokes
	equations consists of the minimization of a cost functional subject to the
	Navier--Stokes equations with the introduction of a control variable
	in the system. The control is supposed to act on the whole spatial domain.
	The formulation of the problem is as follows: given a bounded Lipschitz
	domain $\Omega \subset \mathbb{R}^{d}$, $d\in\{2,3\}$, the
	Navier--Stokes control problem considered in this work is defined as
	\begin{equation}\label{Navier_Stokes_control_cost_functional}
		\min_{\vec{v},\vec{u}}~~J(\vec{v},\vec{u}) =
			\dfrac{1}{2} \int_{\Omega} |\vec{v}(\mathbf{x})-
				\vec{v}_{d}(\mathbf{x})|^2 \: {\rm d}\Omega +
			\dfrac{\beta}{2} \int_{\Omega} |\vec{u}(\mathbf{x})|^2 \:
				{\rm d}\Omega
	\end{equation}
	subject to
	\begin{equation}\label{Navier_Stokes_control_constraints}
		\left\{
			\begin{array}{rl}
				\vspace{1ex}
				- \nu \nabla^2 \vec{v} + \vec{v} \cdot \nabla \vec{v} +
					\nabla p = \vec{u} + \vec{f}(\mathbf{x})
					& \quad \mathrm{in} \; \Omega, \\
				\vspace{1ex}
				-\nabla \cdot \vec{v} (\mathbf{x}) = 0
					& \quad \mathrm{in} \; \Omega, \\
				\vec{v}(\mathbf{x}) = \vec{g}(\mathbf{x})
					& \quad \mathrm{on} \; \partial \Omega,
			\end{array}
		\right.
	\end{equation}
	where $\vec{v}$ and $p$ are the \emph{state velocity} and \emph{state
	pressure}, respectively, and $\vec{u}$ is the \emph{control} (the precise
	spaces in which solutions are sought will be specified later).
	The functions $\vec{f}$ and $\vec{g}$ represent the forcing term acting on
	the system and suitable boundary conditions, respectively. Further,
	$\vec{v}_d$ is the desired (velocity) state, and the parameter
	$\beta>0$ is a regularization parameter. Finally, the parameter
	$\nu>0$ represents the kinematic viscosity of the fluid. The type of
	flow described by the Navier--Stokes equations is influenced by the
	viscosity $\nu$. Specifically, denoting by $L$ and $V$ the characteristic
	length and velocity scales of the flow, respectively, the nature of the
	flow is determined by the Reynolds number, defined as
	$Re = \frac{L \, V }{\nu}$: from small to moderate Reynolds number, the
	flow is laminar, while at high Reynolds number the flow may become
	turbulent. From a point of  view of applications, the distributed
	control of incompressible viscous fluid can be realized, for example, in
	polarizable fluids by applying an electromagnetic force, see
	\cite{Weier_Gerbeth}. Other applications can be  
	found in, e.g., \cite{Gunzburger03} and references therein.

The non-linearity of the PDEs \eqref{Navier_Stokes_control_constraints}
	obviously requires one to employ a non-linear
	iteration in order to obtain an approximation of the solution of a
	Navier--Stokes control problem. The aim of this work is to derive effective
	preconditioners for the linear systems arising from the non-linear iteration
	employed. In \cite{Leveque_Pearson_2022}, Leveque and Pearson considered a
	Picard iteration applied to the control of the (stationary and
	instationary) Navier--Stokes equations. The Picard linearization of the
	first-order optimality conditions of the Navier--Stokes control problem
	results in two coupled \emph{Oseen equations}, to be solved simultaneously.
	By making use of a block-commutator argument, the authors in
	\cite{Leveque_Pearson_2022} were able to derive a robust preconditioner for
	the resulting discretized Oseen equations. The preconditioner showed a mild
	dependence on the viscosity of the fluid, while being robust with respect to
	the mesh-size. With this approach, the authors of
	\cite{Leveque_Pearson_2022} were able
	to solve problems with a rather small viscosity. However, as the viscosity
	decreased,  an increase in the number of Picard
	iterations required for reaching a prescribed reduction on the non-linear
	relative residual was observed. This is not surprising, as
	the contraction factor governing the (linear) convergence of the Picard
	iteration applied to the stationary Navier--Stokes equations
	is known to approach 1 as the Reynolds number approaches its critical
	value $Re_*$ \cite{Kara82,GP83}.
	One way to obtain faster convergence of the non-linear iteration is
	to employ second-order methods, such as Newton's method, which, however,
	necessitates a good initial guess; the smaller the viscosity, the closer
	the initial guess must be to the solution in order for Newton's
	method to converge; the convergence is, eventually, quadratic
	\cite{GP83,Kara82}. While these considerations apply to the forward
	problem, one can expect them to hold also for the control problem.
	
In the context of optimization problems, Newton's method makes use of second
	derivative information and, generally speaking, the linear systems arising
	from Newton's method are more complicated, and harder to solve, than the
	ones arising from Picard's iteration.  For this reason, but also in order
	to avoid computing the full Hessian matrix, use is often made
	of suitable approximations. One class of methods
	obtained in this way are the so called \emph{inexact Newton methods},
	see \cite{Nocedal_Wright}. In this framework,
	one can express the derivatives of the cost functional in terms of the
	Jacobian of the residual. From here, a new approximation of a solution of
	the problem can be obtained by neglecting certain second-order terms arising
	from the Hessian. Additionally, the linear systems involved in the
	computation of the Newton steps need not be solved to high accuracy.

For the inexact Newton iteration applied to the stationary Navier--Stokes
	control problem to be viable, efficient and robust
	preconditioners are needed in the solution of the resulting linear systems.
	Despite the excellent performance of the preconditioner 
	derived in \cite{Leveque_Pearson_2022} for the discretized Picard
	linearization of the first-order optimality conditions, one cannot expect
	this preconditioner to be efficient when applied in the context 
	of an inexact Newton method. In fact, the \emph{Newton matrix} arising in
	a Newton step is not defined on the pressure space, thus it cannot be
	included in the commutator argument. From here, one can expect the
	block-commutator approach to perform poorly for problems with moderate
	to very small viscosity (a fact that is already known for the
	forward problem, see, for instance,
	\cite[Section 9.3.3]{Elman_Silvester_Wathen}).

A remarkably robust approach for solving for the Newton step applied to
	the (forward) Navier--Stokes equations is  augmented Lagrangian-based
	preconditioning. This type of preconditioner was developed in
	\cite{Benzi_Olshanskii_2006} for the solution of the linear systems
	arising from a Picard linearization of the stationary 2D Navier--Stokes
	equations. The crucial ingredient in this approach is a specialized
	multigrid solver used to approximately solve the augmented momentum
	equation. The preconditioner was shown to be robust with respect
	to the mesh-size and, more importantly, to the viscosity of the fluid.
	In \cite{Farrell_Mitchell_Wechsung} this  strategy was extended
	to the solution of the linear systems arising from a Newton
	linearization of the 3D Navier--Stokes equations. Again, 
	the robustness of the preconditioner was confirmed, allowing the
	authors to solve problems with very small viscosity.
	In the present work, we employ an augmented Lagrangian preconditioner
	in the solution of linear systems of saddle-point type arising in the
	inexact Newton iteration applied to the stationary incompressible
	Navier--Stokes control problem. We emphasize that the structure of
	these linear systems is quite different from the one considered in 
	\cite{Benzi_Olshanskii_2006,Farrell_Mitchell_Wechsung}.

The paper is structured as follows. In Section \ref{sec_2}, we introduce the
	non-linear iteration employed in this work for the solution of the
	stationary Navier--Stokes control problem
	\eqref{Navier_Stokes_control_cost_functional}--\eqref{Navier_Stokes_control_constraints}. In Section \ref{sec_3}, we introduce
	the augmented Lagrangian preconditioner employed in this work for solving
	the linear systems arising from the linearization of the problem adopted.
	In Section \ref{sec_4}, we present numerical experiments showing the
	robustness of the proposed methodology. Finally, in Section \ref{sec_5}
	we draw our conclusions and present future developments of this work.

\section{Inexact Newton iteration}\label{sec_2}
In this section, we describe the non-linear iteration that we employ in
	order to solve the stationary incompressible Navier--Stokes control
	problem \eqref{Navier_Stokes_control_cost_functional}--\eqref{Navier_Stokes_control_constraints}.

In order to introduce the non-linear iteration, one has to derive the 
	first-order optimality conditions (also called \emph{Karush--Kuhn--Tucker
	(KKT) conditions}). This can be done in two ways. On one hand, one can first
	derive the optimality conditions in the infinite-dimensional Hilbert
	space where the solutions are sought, then discretize the conditions so
	obtained; this approach is known as \emph{optimize-then-discretize}. On
	the other hand, one can first discretize the cost functional and the
	constraints in \eqref{Navier_Stokes_control_cost_functional}--\eqref{Navier_Stokes_control_constraints},
	then derive the optimality conditions by mean of classical optimization
	theorems in the finite-dimensional setting; this approach is known as
	\emph{discretize-then-optimize}.

In the following, we derive the Newton iteration
	when employing  either the optimize-then-discretize or the
	discretize-then-optimize approach. We opt for this choice
	for two reasons. First, we show the flexibility of the augmented Lagrangian
	preconditioner when either one of the above strategies is employed.
	Second, numerical examples demonstrate the effectiveness and robustness
	of the approach, for a wide range of problem parameters.

The Newton iteration produces
	a sequence of approximations $\vec{v}^{\: (k)}(\mathbf{x}) \approx
	\vec{v}(\mathbf{x})$ and $\vec{u}^{(k)} (\mathbf{x}) \approx
	\vec{u}(\mathbf{x})$ of the state velocity and control
	respectively. The pressures $p(\mathbf{x})$ and $\mu(\mathbf{x})$ arise
	in the optimality conditions as the Lagrange multipliers enforcing the
	adjoint and forward incompressibility constraints. We employ inf--sup stable
	finite element pairs to discretize the velocity-pressure pair. Letting
	$\{\vec{\phi}_i\}_{i=1}^{n_v}$ and $\{\psi_i\}_{i=1}^{n_p}$ be
	the basis functions for an inf--sup stable pair of finite elements for
	velocity and pressure respectively, we introduce the following matrices:
	\begin{align*}
		\mathbf{K} = \left\{ \int_\Omega \nabla \vec{\phi}_i : \nabla
			\vec{\phi}_j \, {\rm d}\Omega \right\}_{i,j=1}^{n_v}, & \qquad
		\mathbf{N}^{(k)} = \left\{ \int_\Omega ( \vec{v}^{\: (k)} \cdot \nabla
			\vec{\phi}_j ) \cdot \vec{\phi}_i \,
			{\rm d}\Omega \right\}_{i,j=1}^{n_v}, \\
		\mathbf{M} = \left\{ \int_\Omega \vec{\phi}_i \cdot
			\vec{\phi}_j \, {\rm d}\Omega \right\}_{i,j=1}^{n_v}, & \qquad
		B = \left\{ - \int_\Omega \psi_i  \nabla \cdot
			\vec{\phi}_j \, {\rm d}\Omega
			\right\}_{i = 1, \ldots, n_p}^{j = 1, \ldots, n_v}.
	\end{align*}
	The matrix $\mathbf{K}$ is generally referred to as a
	(vector-)stiffness matrix, and the matrix $\mathbf{M}$ is referred
	to as a (vector-)mass matrix; both of these matrices are symmetric
	positive definite (s.p.d.). The matrix $\mathbf{N}^{(k)}$ is referred to
	as a (vector-)convection matrix, and is skew-symmetric (i.e.
	$\mathbf{N}^{(k)}+(\mathbf{N}^{(k)})^\top=0$) if the ``wind''
	$\vec{v}^{\: (k)}$ is incompressible; finally, the matrix $B$ is
	the (negative) divergence matrix. In addition to the previous
	matrices, we also introduce the following:
	\begin{equation}\label{Newton_derivative_matrix}
		\mathbf{H}^{(k)} = \left\{ \int_\Omega (\vec{\phi}_j \cdot \nabla
			\vec{v}^{\: (k)}) \cdot \vec{\phi}_i \, {\rm d}\Omega
			\right\}_{i,j=1}^{n_v},
	\end{equation}
	which is the matrix arising in Newton's method from linearizing
	the convection term.

At high Reynolds number, the problem is convection-dominated. This requires
	us to make use of a stabilization procedure. In the following,
	$\mathbf{W}^{(k)}$ denotes a possible stabilization matrix for the
	convection operator. It should be mentioned that the matrix
	$\mathbf{W}^{(k)}$ represents a differential operator that is not related
	to the physics, and is introduced solely to enhance coercivity (that is,
	increase the positivity of the real part of the eigenvalues) of the
	discretization, thereby allowing it to be stable. Following the
	discussions in \cite{Leveque_Pearson_2021, Pearson_Wathen_2013}, in some
	of the numerical experiments presented in this work we employ the
	Local Projection Stabilization (LPS) approach described in
	\cite{Becker_Braack, Becker_Vexler, Braack_Burman}. For an analysis of
	the order of convergence of the LPS applied to the forward Oseen problem,
	we refer to \cite{Matthies_Tobiska}. For other possible stabilizations
	applied to the forward Oseen problem, see, for instance,
	\cite{Brooks_Hughes, Franca_Frey, Johnson_Saranen, Tobiska_Lube}.

We would like to mention here that the choice of employing the
	LPS stabilization has been made for an important reason. When solving the
	control of the convection--diffusion equation, the discretize-then-optimize and
	the optmize-then-discretize approaches do not commute if a different
	stabilization is employed, see \cite{Collis_Heinkenschloss}. However, as shown
	in \cite{Becker_Vexler}, when employing the LPS approach to the control of the
	stationary convection--diffusion equation, the discretization and optimization
	steps commute. Since at each Newton iteration one has to solve a
	convection--diffusion control problem, the choice of the LPS as stabilization
	is very natural.

The LPS formulation is defined as follows. Given $\pi_h$, an $L^2$-orthogonal
	(discontinuous) projection operator onto the finite dimensional
	space $\overline{V}$ defined on patches of $\Omega$, where
	by a patch we mean a union of cells sharing a vertex,
	we consider the \emph{fluctuation operator} $\kappa_h=\text{Id} -\pi_h$.
	Here, $\text{Id}$ denotes the identity operator. Then, the matrix
	$\mathbf{W}^{(k)}$ is defined as
	\begin{equation}\label{Navier_Stokes_LPS_stab_matrix}
		\mathbf{W}^{(k)} = \left\{ \delta^{(k)} \int_\Omega
			\kappa_h(\vec{v}^{\: (k)} \cdot \nabla \vec{\phi}_i) \cdot
			\kappa_h(\vec{v}^{\: (k)} \cdot \nabla \vec{\phi}_j)
			\, {\rm d}\Omega \right\}_{i,j=1}^{n_v},
	\end{equation}
	where $\delta^{(k)} \ge 0$ is a stabilization parameter. As in
	\cite{Becker_Braack}, we define the projection by means of
	\begin{displaymath}
		\left.\pi_h(q)\right|_{\mathtt{P}_m} = \dfrac{1}{|\mathtt{P}_m|}
			\int_{\mathtt{P}_m} q \: \mathrm{d} \mathtt{P}_m, \quad
			\forall q \in L^2(\Omega),
	\end{displaymath}
	where $|\mathtt{P}_m|$ is the (Lebesgue) measure of the patch
	$\mathtt{P}_m$. We refer the reader to \cite{Matthies_Tobiska} for
	the motivation of the choice of this projection and for the possible
	choices of the finite dimensional space $\overline{V}$. Finally, as
	in \cite[p.\,253]{Elman_Silvester_Wathen}
	the stabilization parameter is taken to be
	\begin{displaymath}
		\delta^{(k)}_m =
			\left\{
				\begin{array}{lcl}
					\vspace{1ex}
					\dfrac{h_m}{2 \| \vec{v}^{\: (k)}_m \| }
						\left(1 - \dfrac{1}{Pe_m}\right) &
						\quad \mathrm{if} \: Pe_m > 1, \\
					0 & \quad \mathrm{if} \: Pe_m \leq 1,
				\end{array}
			\right.
	\end{displaymath}
	where $\| \vec{v}^{\: (k)}_m \|$ is the Euclidean norm of $\vec{v}^{\: (k)}$
	at the patch centroid, $h_m$ is a measure of the patch length in the
	direction of the wind, and $Pe_m=\| \vec{v}_m^{\: (k)} \| h_m / (2 \nu)$
	is the patch P\'{e}clet number. For the discretization to be stable,
	we expect $Pe_m$ to be less than 1. We would like to mention that the choice
	we made of the stabilization parameter is only a heuristic.

In the following, we set $\mathbf{D}^{(k)} = \nu \mathbf{K} + \mathbf{N}^{(k)}
	+ \mathbf{W}^{(k)}$. Further, we define
	$$V : = \{\vec{v} \in
	\mathcal{H}^1(\Omega)^d \: | \: \vec{v} = \vec{g} \; \mathrm{on} \;
	\partial \Omega \}, \quad
	V_0 : = \{\vec{v} \in \mathcal{H}^1(\Omega)^d \: |
	\: \vec{v} = \vec{0} \; \mathrm{on} \; \partial \Omega \},\quad  
	Q:= L^2_0(\Omega),$$
	 with $\mathcal{H}^1(\Omega)^d$ the usual Sobolev space
	of square-integrable functions in $\mathbb{R}^d$ with square-integrable
	weak derivatives.

\subsection{Optimize-then-Discretize}\label{sec_2_1}
We begin by considering the optimize-then-discretize approach, which
	makes use of the so called formal Lagrangian method, see
	\cite[Section 2.10]{Troltzsch}.

First, we derive the continuous KKT conditions. Introducing the
	adjoint velocity $\vec{\zeta}$ and the adjoint
	pressure $\mu$, one may consider the Lagrangian associated with
	\eqref{Navier_Stokes_control_cost_functional}--\eqref{Navier_Stokes_control_constraints},
	as in \cite{Posta_Roubicek}. Then, one can write the KKT conditions as
	\begin{equation}\label{Navier_Stokes_control_KKT_conditions}
		\left\{
			\hspace{-2ex}
			\begin{array}{rl}
				\vspace{0.25ex}
				\left.
					\begin{array}{rl}
						- \nu\nabla^2 \vec{\zeta} - \vec{v} \cdot \nabla
							\vec{\zeta} + (\: \nabla \vec{v} \:) ^\top
							\vec{\zeta} + \nabla \mu \: = \: \vec{v}_d -
							\vec{v} & 
							\mathrm{in} \; \Omega \\
						- \nabla \cdot \vec{\zeta}(x) = 0 &
							\mathrm{in} \; \Omega \\
						\vec{\zeta}(x) = \vec{0} & 
							\mathrm{on} \; \partial \Omega
					\end{array}
				\right\} & 
				\left.
					\hspace{-1.5ex}
					\begin{array}{c}
						\vspace{0.5ex}
						\mathrm{adjoint}\\
						\mathrm{equations}
					\end{array}
				\right.\\
					\begin{array}{rl}
						\beta \vec{u} - \vec{\zeta} = 0 & 
							\mathrm{in} \; \Omega \phantom{\:\:\:\:\:\:\:\,}
					\end{array}
				 & 
				\left.
					\hspace{-1.5ex}
					\begin{array}{c}
						\vspace{0.5ex}
						\mathrm{gradient}\\
						\mathrm{equation}
					\end{array}
				\right.\\
				\left.
					\begin{array}{rl}
						- \nu\nabla^2 \vec{v} + \vec{v} \cdot \nabla \vec{v}
							+ \nabla p = \vec{u} +
							\vec{f} &
						\mathrm{in} \; \Omega \\
						- \nabla \cdot \vec{v}(x) = 0 & 
							\mathrm{in} \; \Omega \\
						\vec{v}(x) = \vec{g}(x) & 
							\mathrm{on} \; \partial \Omega
					\end{array}
				\right\} & 
				\left.
					\hspace{-1.5ex}
					\begin{array}{c}
						\vspace{0.5ex}
						\mathrm{state}\\
						\mathrm{equations}
					\end{array}
				\right.
			\end{array}
		\right.
	\end{equation}
	The expression of the gradient equation motivates us to take $\vec{v}$
	(thus, $\vec{\zeta}$) and $\vec{u}$ in the same space, so that we can
	eliminate $\vec{u}$ from the system.

Problem \eqref{Navier_Stokes_control_KKT_conditions} consists of a coupled
	system of non-linear PDEs. We employ an inexact Newton method in order
	to derive an approximation of the solutions. The goal is to write the
	Newton step for solving the KKT conditions in
	\eqref{Navier_Stokes_control_KKT_conditions}, and then neglect the
	\emph{curvature term} of the full Hessian. In order to write the
	Newton system, we follow the work in \cite{Hinze_Kunisch}.

Given $\vec{v}^{\:( k)}\in V$, $p^{(k)} \in Q$, $\vec{\zeta}^{\: (k)} \in V_0$,
	$\mu^{(k)} \in Q$  the current approximations to $\vec{v}$, $p$,
	$\vec{\zeta}$, and $\mu$, respectively, one can write the Newton
	iterates as
	\begin{displaymath}
		\begin{array}{c}
			\vec{v}^{\: (k+1)}=\vec{v}^{\: (k)}+ \vec{\delta v}^{\: (k)},
				\qquad p^{(k+1)} = p^{(k)} + \delta p ^{(k)}, \\
			\vec{\zeta}^{\: (k+1)}= \vec{\zeta}^{\: (k)} +
				\vec{\delta \zeta}^{\: (k)}, \qquad
			\mu^{(k+1)}=\mu^{(k)} + \delta \mu^{(k)},
		\end{array}
	\end{displaymath}
	where the Newton corrections are solutions of the following system
	of PDEs (written in weak form):
	\begin{equation}\label{Opt_then_Discr_Newton_step}
		\left\{
			\begin{array}{rl}
				\nu(\nabla \vec{\delta \zeta}^{\: (k)},
					\nabla \vec{w}_1) 
					-  (\vec{v}^{\: (k)} \cdot \nabla \vec{\delta
					\zeta}^{\: (k)}, \vec{w}_1)
					+ ((\: \nabla \vec{v}^{\: (k)} \:) ^\top
					\vec{\delta \zeta}^{\: (k)}, \vec{w}_1) \\
				- ( \delta \mu^{(k)}, \nabla \cdot \vec{w}_1) \! + \!
					(\vec{\delta v}^{\: (k)},\vec{w}_1)
					+ (\vec{\zeta}^{\: (k)} \cdot \nabla
						\vec{\delta v}^{\: (k)},
						\vec{w}_1) \\
					+ (\vec{\delta v}^{\: (k)} \cdot \nabla
					\vec{\zeta}^{\: (k)}, \vec{w}_1) \!  &
					\hspace{-0.6em} = \vec{R}_1^{\: (k)}, \\
				- (q_2,\nabla \cdot \vec{\delta \zeta}^{\: (k)})
					& \hspace{-0.6em} = r_2^{(k)},\\
				\nu (\nabla \vec{\delta v}^{\: (k)}, \nabla \vec{w}_2 ) \!
					+ ( \vec{v}^{\: (k)} \cdot \nabla \vec{\delta v}^{\: (k)},
					\vec{w}_2) \!
					+ ( \vec{\delta v}^{\: (k)} \cdot \nabla \vec{v}^{\: (k)},
					\vec{w}_2) \\
					- ( \delta p^{(k)}, \nabla  \cdot \vec{w}_2 ) -
					\frac{1}{\beta} ( \vec{\delta \zeta}^{\: (k)},\vec{w}_2)
					&
					\hspace{-0.6em} = \vec{R}_2^{\: (k)}, \\
				- (q_1,\nabla \cdot \vec{\delta v}^{\: (k)}) &
					\hspace{-0.6em} = r_1^{(k)},
			\end{array}
		\right.
	\end{equation}
	for any $(\vec{w}_1, \vec{w}_2, q_1, q_2) \in V_0 \times V_0 \times Q
	\times Q$. The residuals $\vec{R}_1^{\: (k)}$,
	$r_1^{(k)}$, $\vec{R}_2^{\: (k)}$, $r_2^{(k)}$ are given by
	\begin{equation}\label{Opt_then_Discr_residual}
		\!\!
		\left\{
			\!\!
			\begin{array}{rl}
				\vspace{0.125ex}
				\vec{R}_1^{\: (k)} \! = & \hspace{-0.8em} (\vec{v}_d, \vec{w}_1)
					- (\vec{v}^{\: (k)},\vec{w}_1) - \nu (\nabla
					\vec{\zeta}^{\: (k)},\nabla \vec{w}_1) + (\vec{v}^{\: (k)}
					\cdot \nabla \vec{\zeta}^{\: (k)},\vec{w}_1) \\
				\vspace{0.125ex}
					& - ((\: \nabla \vec{v}^{\: (k)} \:) ^\top
					\vec{\zeta}^{\: (k)}, \vec{w}_1) + ( \mu^{(k)}, \nabla
					\cdot \vec{w}_1), \\
				r_2^{(k)} \! = & \hspace{-0.8em} (q_2, \nabla \cdot
					\vec{\zeta}^{\:( k)}), \\
				\vec{R}_2^{\: (k)} \! = & \hspace{-0.8em} (\vec{f},\vec{w}_2) \!
					- \! \nu(\nabla \vec{v}^{\: (k)},\nabla \vec{w}_2) \! - \!
					(\vec{v}^{\: (k)} \cdot \nabla \vec{v}^{\: (k)}, \vec{w}_2)
					\! + \! ( p^{(k)}, \nabla \cdot \vec{w}_2) \!\\
					& + \frac{1}{\beta}( \vec{ \zeta}^{\: (k)},\vec{w}_2), \\
				\vspace{0.125ex}
				r_1^{(k)} \! = & \hspace{-0.8em} (q_1,\nabla \cdot
					\vec{v}^{\: (k)}).
			\end{array}
		\right.
	\end{equation}

We have now to derive the discrete system to solve at each non-linear iteration.
	Before doing this, we derive an inexact Newton method by neglecting the
	coupling terms $(\vec{\zeta}^{\: (k)} \cdot \nabla
	\vec{\delta v}^{\: (k)}, \vec{w}_1)$ and $(\vec{\delta v}^{\: (k)}
	\cdot \nabla \vec{\zeta}^{\: (k)}, \vec{w}_1)$ in
	\eqref{Opt_then_Discr_Newton_step}, as done in \cite{Benzi_Haber_Taralli},
	for instance. As we will show below, dropping these terms allows us to
	derive effective and robust preconditioners for the linearized
	inexact Newton system, while still maintaining second-order convergence
	of the method. The inexact Newton iteration we consider is given by the
	following system of PDEs:
	\begin{equation}\label{Opt_then_Discr_Inexact_Newton_step}
		\left\{
			\begin{array}{rl}
				(\vec{\delta v}^{\: (k)},\vec{w}_1)
					+ \nu(\nabla \vec{\delta \zeta}^{\: (k)},
					\nabla \vec{w}_1)
					-  (\vec{v}^{\: (k)} \cdot \nabla \vec{\delta
					\zeta}^{\: (k)}, \vec{w}_1) \\
					+ ((\: \nabla \vec{v}^{\: (k)} \:) ^\top
					\vec{\delta \zeta}^{\: (k)}, \vec{w}_1)
					- ( \delta \mu^{(k)}, \nabla \cdot \vec{w}_1) &
					\hspace{-0.6em} = \vec{R}_1^{\: (k)}, \\
				- (q_2,\nabla \cdot \vec{\delta \zeta}^{\: (k)})
					& \hspace{-0.6em} = r_2^{(k)}, \\
				\nu (\nabla \vec{\delta v}^{\: (k)}, \nabla \vec{w}_2 ) \!
					+ ( \vec{v}^{\: (k)} \cdot \nabla \vec{\delta v}^{\: (k)},
					\vec{w}_2) \!
					+ ( \vec{\delta v}^{\: (k)} \cdot \nabla \vec{v}^{\: (k)},
					\vec{w}_2) \\
					- ( \delta p^{(k)}, \nabla  \cdot \vec{w}_2 ) -
					\frac{1}{\beta} ( \vec{\delta \zeta}^{\: (k)},\vec{w}_2)
					&
					\hspace{-0.6em} = \vec{R}_2^{\: (k)}, \\
				- (q_1,\nabla \cdot \vec{\delta v}^{\: (k)}) &
					\hspace{-0.6em} = r_1^{(k)},
			\end{array}
		\right.
	\end{equation}
	for any $(\vec{w}_1, \vec{w}_2, q_1, q_2) \in V_0 \times V_0 \times Q
	\times Q$.

Next, assuming that an inf-sup stable finite element discretization is 
applied to  \eqref{Opt_then_Discr_Inexact_Newton_step}, we 
can write the (finite-dimensional) inexact Newton step in the form
	\begin{displaymath}
		\begin{array}{c}
			\bm{v}^{\: (k+1)}=\bm{v}^{\: (k)}+ \bm{\delta v}^{\: (k)},
				\qquad \bm{p}^{(k+1)} = \bm{p}^{(k)} + \bm{\delta p }^{(k)}, \\
			\bm{\zeta}^{\: (k+1)}= \bm{\zeta}^{\:( k)} +
				\bm{\delta \zeta}^{\:( k)}, \qquad
				\bm{\mu}^{(k+1)}=\bm{\mu}^{(k)} + \bm{\delta \mu}^{(k)},
		\end{array}
	\end{displaymath}
	where the inexact Newton corrections satisfy the following linear
	algebraic system:
	\begin{equation}\label{Opt_then_Discr_discrete_Inexact_Newton_system}
		\underbrace{
			\left[
				\begin{array}{cc}
				{\mathcal F}^{(k)}_{\mathrm{OD}} & {\mathcal B}^\top \\
				         {\mathcal B}  &  {\mathcal O}
				\end{array}
			\right]
		}_{\mathcal{A}^{(k)}_{\mathrm{OD}}}
		\left[
			\begin{array}{c}
				\bm{\delta v}^{(k)}\\
				\bm{\delta \zeta}^{(k)}\\
				\bm{\delta \mu}^{(k)}\\
				\bm{\delta p}^{(k)}
			\end{array}
		\right]=
		\left[
			\begin{array}{c}
				\bm{R}_1^{(k)}\\
				\bm{R}_2^{(k)}\\
				\bm{r}_1^{(k)}\\
				\bm{r}_2^{(k)}
			\end{array}
		\right],
	\end{equation}
	where, by setting $\mathbf{M}_\beta = \frac{1}{\beta}\mathbf{M}$
	and $\mathbf{D}^{(k)}_\mathrm{adj} = \nu \mathbf{K}
	- \mathbf{N}^{(k)} + \mathbf{W}^{(k)}$, we have
	\begin{equation}\label{Phi_Psi_Theta_Opt_then_Discr}
		{\mathcal F}^{(k)}_\mathrm{OD} =
		\left[ \!\!
			\begin{array}{cc}
				\mathbf{M} &
					\mathbf{D}^{(k)}_{\mathrm{adj}} + (\mathbf{H}^{(k)})^\top\\
				\mathbf{D}^{(k)} + \mathbf{H}^{(k)} & -\mathbf{M}_\beta
			\end{array}
		\!\! \right], \quad
		{\mathcal B} =
		\left[ \!\!
			\begin{array}{cc}
				B & 0\\
				0 & B
			\end{array}
		\!\! \right], \quad
		{\mathcal O} =
		\left[ \!\!
			\begin{array}{cc}
				0 & 0\\
				0 & 0
			\end{array}
		\!\! \right].
	\end{equation}
	The discrete residuals
	in \eqref{Opt_then_Discr_discrete_Inexact_Newton_system} are given by
	\begin{equation}\label{Opt_then_Discr_discrete_Residual}
		\left\{
			\begin{array}{rll}
				\bm{R}_1^{\: (k)} = & \hspace{-0.6em} \mathbf{M} \: \bm{v}_d -
					\mathbf{M} \: \bm{v}^{\: (k)} -
					\mathbf{D}^{(k)}_\mathrm{adj} \: \bm{\zeta}^{\: (k)}
					- B^\top \bm{\mu}^{(k)}  - \bm{\omega}^{(k)},\\
				\bm{R}_2^{\: (k)} = & \hspace{-0.6em} \bm{f} -
					\mathbf{D}^{(k)} \: \bm{v}^{\: (k)} - B^\top \:
					\bm{p}^{(k)} + \mathbf{M}_\beta\: \bm{ \zeta}^{\: (k)}, \\
				\bm{r}_1^{(k)} = & \hspace{-0.6em} - B \bm{v}^{\: (k)}, \\
				\bm{r}_2^{(k)} = & \hspace{-0.6em} - B \bm{\zeta}^{\: (k)},
			\end{array}
		\right.
	\end{equation}
	where $\bm{v}_d$ is the vector
	corresponding to the discretized desired state $\vec{v}_d$, $\bm{f}$ is
	the vector corresponding to the discretized force function $f$, and
	$\bm{\omega}^{(k)}=\{ \big( (\nabla \vec{v}^{\:(k)} )^\top 
	\vec{\zeta}^{\: (k)}, \vec{\phi}_i \big) \}_{i=1}^{n_v}$.

\subsection{Discretize-then-Optimize}\label{sec_2_2}
We now present the inexact Newton iteration for the
	discretize-then-optimize approach. We suppose that we have approximations
	$\vec{v}^{\:(k)}$ and $\vec{\zeta}^{\: (k)}$ of the state velocity
	$\vec{v}$ and of the adjoint velocity $\vec{\zeta}$, respectively.

We first discretize the cost functional in
	\eqref{Navier_Stokes_control_cost_functional} and the constraints in
	\eqref{Navier_Stokes_control_constraints}. In order to derive a second-order
	method, we need to employ a Newton discretization of the convection
	operator. Specifically, we have
	\begin{displaymath}
		\min_{\bm{v},\bm{u}}~~\mathbf{J}(\bm{v},\bm{u}) =
			\frac{1}{2} ( \bm{v} - \bm{v}_d )^\top \mathbf{M} ( \bm{v} -
				\bm{v}_d ) + \frac{\beta}{2} \bm{u}^\top \mathbf{M} \bm{u},
	\end{displaymath}
	subject to
	\begin{displaymath}
		\left\{
			\begin{array}{r}
				\vspace{1ex}
				(\mathbf{D}^{(k)} + \mathbf{H}^{(k)})\bm{v} + B^\top \bm{p}
					- \mathbf{M} \bm{u} = \bm{f}, \\
				\vspace{1ex}
				B \bm{v} = \bm{0}, \\
				\bm{v}_{\mathrm{b}} = \bm{g},
			\end{array}
		\right.
	\end{displaymath}
	where $\bm{v}_{\mathrm{b}}$ are the components of the vector $\bm{v}$
	related to the boundary nodes. By employing classical constrained
	optimization theory, one can derive the first-order optimality conditions
	for the problem above. Specifically, the KKT conditions are given by
	the following system of equations:
	\begin{displaymath}
		\left\{
			\begin{array}{rl}
				(\mathbf{D}^{(k)} + \mathbf{H}^{(k)}) \bm{v} + B^\top \bm{p} -
					\frac{1}{\beta} \mathbf{M} \bm{\zeta} &  \!\! = \bm{f}, \\
				\vspace{1ex}
				B \bm{v} & \!\! = \bm{0}, \\
				\mathbf{M} \bm{v} + (\mathbf{D}^{(k)} + \mathbf{H}^{(k)})^\top
					\bm{\zeta} + B^\top \bm{\mu}
					& \!\! =\mathbf{M} \bm{v}_d, \\
				B \bm{\zeta} & \!\! = \bm{0},
			\end{array}
		\right.
	\end{displaymath}
	where we have substituted the gradient equation
	$\beta \bm{u} = \bm{\zeta}$ into the state equation.

As a last step, we derive the inexact Newton iteration. This is given by
	\begin{displaymath}
		\begin{array}{c}
			\bm{v}^{\: (k+1)}=\bm{v}^{\: (k)}+ \bm{\delta v}^{\: (k)},
				\qquad \bm{p}^{(k+1)} = \bm{p}^{(k)} + \bm{\delta p }^{(k)}, \\
			\bm{\zeta}^{\: (k+1)}= \bm{\zeta}^{\:( k)} +
				\bm{\delta \zeta}^{\:( k)}, \qquad
				\bm{\mu}^{(k+1)}=\bm{\mu}^{(k)} + \bm{\delta \mu}^{(k)},
		\end{array}
	\end{displaymath}
	where the inexact Newton corrections satisfy the following linear system:
	\begin{equation}\label{Discr_then_Opt_discrete_Inexact_Newton_system}
		\underbrace{
			\left[
				\begin{array}{cc}
					\mathcal{F}^{(k)}_{\mathrm{DO}} & {\mathcal B}^\top\\
					{\mathcal B} & {\mathcal O}
				\end{array}
			\right]
		}_{\mathcal{A}^{(k)}_{\mathrm{DO}}}
		\left[
			\begin{array}{c}
				\bm{\delta v}^{(k)}\\
				\bm{\delta \zeta}^{(k)}\\
				\bm{\delta \mu}^{(k)}\\
				\bm{\delta p}^{(k)}
			\end{array}
		\right]=
		\left[
			\begin{array}{c}
				\bm{R}_1^{(k)}\\
				\bm{R}_2^{(k)}\\
				\bm{r}_1^{(k)}\\
				\bm{r}_2^{(k)}
			\end{array}
		\right],
	\end{equation}
	with the discrete non-linear residual given as in
	\eqref{Opt_then_Discr_discrete_Residual}, the matrices $\mathcal B$ and
	$\mathcal O$ defined as in \eqref{Phi_Psi_Theta_Opt_then_Discr}, and
	\begin{equation}\label{Phi_Discr_then_Opt_Inexact_Newton}
		\mathcal{F}^{(k)}_{\mathrm{DO}} =
		\left[
			\begin{array}{cc}
				\mathbf{M} &
					(\mathbf{D}^{(k)} + \mathbf{H}^{(k)})^\top\\
				\mathbf{D}^{(k)} + \mathbf{H}^{(k)} & -\mathbf{M}_\beta
			\end{array}
		\right].
	\end{equation}

Note that the linear system given in
	\eqref{Discr_then_Opt_discrete_Inexact_Newton_system} does not correspond
	to the full Newton iteration, as in the $(1,1)$-block we do not
	include second-order information. In fact, introducing the function
	\begin{displaymath}
		\mathcal{G}(\bm{v}, \bm{\zeta}) = \bigintsss_{\Omega}
			\left(\sum_{i=1}^{n_v} \bm{v}_i \vec{\phi}_i \right)
			\cdot \left(\sum_{i=1}^{n_v} \bm{v}_i \nabla \vec{\phi}_i \right)
			\cdot \left(\sum_{i=1}^{n_v} \bm{\zeta}_i \vec{\phi}_i \right)
			\mathrm{d} \Omega,
	\end{displaymath}
	the full Newton system would be given by
	\eqref{Discr_then_Opt_discrete_Inexact_Newton_system} with the $(1,1)$-block
	of the matrix $\mathcal{F}^{(k)}_\mathrm{DO}$ given by $\mathbf{M}
	+ \nabla_{\bm{v}\bm{v}} \mathcal{G}(\bm{v}^{(k)}, \bm{\zeta}^{(k)})$.
	After some manipulations, one can write 
	\begin{displaymath}
		\nabla_{\bm{v}\bm{v}} \mathcal{G}(\bm{v}^{(k)}, \bm{\zeta}^{(k)}) =
			\mathbf{N}^{(k)}_{\bm{\zeta}} + \mathbf{H}^{(k)}_{\bm{\zeta}},
	\end{displaymath}
	with
	\begin{align*}
		\mathbf{N}^{(k)}_{\bm{\zeta}} = \left\{ \int_\Omega (
			\bm{\zeta}^{\: (k)} \cdot \nabla
			\vec{\phi}_j ) \cdot \vec{\phi}_i \, {\rm d}\Omega
			\right\}_{i,j=1}^{n_v},
		& \qquad
		\mathbf{H}^{(k)}_{\bm{\zeta}} = \left\{ \int_\Omega (\vec{\phi}_j
			\cdot \nabla \bm{\zeta}^{\: (k)}) \cdot \vec{\phi}_i
			\, {\rm d}\Omega \right\}_{i,j=1}^{n_v}.
	\end{align*}

\section{Preconditioning approach}\label{sec_3}
In this section, we present the preconditioning strategy that we adopt for
	solving the sequence of linear systems arising in the non-linear iterations.

We are interested in the solution of linear systems of the form
	\begin{equation}\label{generalized_saddle_point}
		\underbrace{\left[
			\begin{array}{cc}
				\Phi & \Psi_1^\top\\
				\Psi_2 & - \Theta
			\end{array}
		\right]}_{\mathcal{A}}
		\left[
			\begin{array}{c}
				\bm{y}_1\\
				\bm{y}_2
			\end{array}
		\right]
		=
		\left[
			\begin{array}{c}
				\bm{b}_1\\
				\bm{b}_2
				\end{array}
		\right],
	\end{equation}
	with invertible $\Phi$ and a possibly non-zero $\Theta$. Systems of
	this form are referred to as generalized saddle-point problems
	in \cite{Benzi_Golub_Liesen}. Note that, for the problems we are interested
	in, both the whole linear system (see
	\eqref{Opt_then_Discr_discrete_Inexact_Newton_system} and
	\eqref{Discr_then_Opt_discrete_Inexact_Newton_system})
	and the subsystems associated with its $(1,1)$-block 
	(see \eqref{Phi_Psi_Theta_Opt_then_Discr} and
	\eqref{Phi_Discr_then_Opt_Inexact_Newton}) are generalized
	saddle-point problems. We also observe that the discretize-then-optimize
	strategy leads to symmetric systems ($\Phi = \Phi^\top$ and
	$\Psi_1 = \Psi_2$), while the optimize-then-discretize one does not since,
	generally speaking,
	$\mathbf{D}^{(k)}_{\mathrm{adj}}\ne (\mathbf{D}^{(k)})^\top$, e.g.~due to lack of exact incompressibility of the state velocity approximation.

As is well-known, an ``ideal" optimal preconditioner for the system in
	\eqref{generalized_saddle_point} is given by either of the following
	block-triangular matrices:
	\begin{equation}\label{optimal_precond}
		\mathcal{P}_1 = \left[
			\begin{array}{cc}
				\Phi & 0\\
				\Psi_2 & -S
			\end{array}
		\right],
		\qquad
		\mathcal{P}_2 = \left[
			\begin{array}{cc}
				\Phi & \Psi_1^\top\\
				0 & -S
			\end{array}
		\right],
	\end{equation}
	where the $(1,1)$-block of the preconditioner $\mathcal{P}_i$,
	$i=1,2$, is the $(1,1)$-block of the matrix $\mathcal{A}$ considered, with the
	$(2,2)$-block $S= \Theta + \Psi_2 \Phi^{-1} \Psi_1^\top$ being the 
	(negative) Schur complement. The optimality of the preconditioners
	derives from the spectral properties of the preconditioned matrices.
	Indeed, assuming that the Schur complement $S$ is invertible,
	it can be proved that $\lambda(\mathcal{P}_i^{-1} \mathcal{A}) = 
	\left\{ 1 \right\}$, for $i=1,2$, with minimal polynomial of degree two;
	see, for instance, \cite{Ipsen01, Murphy_Golub_Wathen}. Therefore, a
	suitable minimum residual method would converge in at most two iterations,
	in exact arithmetic.

In practical applications, the linear system in \eqref{generalized_saddle_point}
	can be very large and $\Phi^{-1}$ full, thus, even forming the Schur
	complement $S$ is often unfeasible. For this reason, many of the most
	effective preconditioners are based on approximations of $\Phi$ and $S$.
	Specifically, given invertible approximations $\widetilde{\Phi}$ and
	$\widetilde{S}$ of $\Phi$ and $S$, respectively, rather than considering
	the ideal preconditioners $\mathcal{P}_1$ and $\mathcal{P}_2$,
	one instead employs approximations of the form
	\begin{equation}\label{approx_preconditioners}
		\widetilde{\mathcal{P}}_1 = \left[
			\begin{array}{cc}
				\vspace{1.05ex}
				\widetilde{\Phi} & 0\\
				\Psi_2 & -\widetilde{S}
			\end{array}
		\right],
		\qquad
		\widetilde{\mathcal{P}}_2 = \left[
			\begin{array}{cc}
				\vspace{1.05ex}
				\widetilde{\Phi} & \Psi_1^\top\\
				0 & -\widetilde{S}
			\end{array}
		\right].
	\end{equation}
	Obviously, the approximations $\widetilde{\Phi}$ and
	$\widetilde{S}$ are problem-dependent. In some cases they are defined
	implicitly through some inexact solution procedure for linear systems
	involving the matrices $\Phi$ and $S$.

As mentioned, for incompressible Navier--Stokes control problems the
	$(1,1)$-block $\Phi$ of the linear system arising in the non-linear
	iteration is itself a (generalized) saddle-point system. Linear systems with
	a similar ``nested" structure arise frequently in the context of
	PDE-constrained optimization problems, and many researchers
	have devoted their efforts in devising robust solvers for such
	linear systems. We refer the interested reader to
	\cite{Axelsson_Farouq_Neytcheva, Heidel_Wathen,
	Qiu_vanGijzen_vanWingerden_Verhaegen_Vuik, Rees_Dollar_Wathen,
	Schoberl_Zulehner, Stoll_Breiten, Zulehner} and the references therein.
	
The most delicate task in the derivation of the preconditioners 
	\eqref{approx_preconditioners} for the (sequence of) linear system(s)
	\eqref{Opt_then_Discr_discrete_Inexact_Newton_system} and
	\eqref{Discr_then_Opt_discrete_Inexact_Newton_system} is finding
	an approximation for the Schur complement $S$. In
	\cite{Leveque_Pearson_2022}, the authors have derived
	a (heuristic) approximation of the Schur complement of the linear systems
	arising in the Picard iteration for the solution of (stationary and
	instationary) Navier--Stokes control problems, based on a commutator
	argument. This strategy has been shown to be very effective and robust for a
	wide range of parameters and problems, see \cite{Danieli_Southworth_Wathen,
	Leveque_Bergamaschi_Martinez_Pearson, Leveque_Pearson_2022}. However, as
	we will see below, when applied to a Newton iteration for the Navier--Stokes
	control problem, the commutator argument presents some limitations. We also
	mention here the work \cite{Xu_Wang}, where the authors derived a rotated
	block-diagonal preconditioner for the solution of the Picard linearization
	of the stationary Navier--Stokes control problem, when employing a
	discretize-then-optimize approach.

In the following sections, we describe the strategies employed for approximating
	the main blocks of the preconditioners we use. In what follows, we
	denote with $I_{n} \in \mathbb{R}^{n \times n}$ the identity matrix.
	Further, $K_p$ and $M_p$ denote stiffness and mass matrices on
	the pressure space, respectively.

\subsection{Matching strategy}\label{sec_3_1}
We first introduce the approximation of the $(1,1)$-block of
	the (whole) linear systems arising from the non-linear iteration.
	The preconditioner we adopt employs the \emph{matching strategy}
	\cite{Pearson_Wathen_2012} for approximating the (inner) Schur complement
	of the $(1,1)$-block, which couples the state and adjoint velocities
	and is of saddle-point type, as in \eqref{generalized_saddle_point}.
	The strategy employed for approximating
	the outer Schur complement (the one of the whole system, where also the
	incompressibility constraints are considered) will be discussed below.
	Although the matching strategy was devised for symmetric indefinite systems,
	in our derivation we consider a matrix $\mathcal{A}$ of the form given in
	\eqref{generalized_saddle_point}, where the $(1,2)$-block may not be
	the transpose of the $(2,1)$-block. This is because the linear
	system \eqref{Phi_Psi_Theta_Opt_then_Discr} arising from the
	optimize-then-discretize strategy, as already observed, may not be
	symmetric. Further, we will assume that the $(1,1)$-block $\Phi$ is s.p.d.,
	and that the $(2,2)$-block is a positive multiple of the $(1,1)$-block,
	namely $\Theta=\frac{1}{\beta}\Phi$, for some $\beta>0$.

As we discussed above, an ideal preconditioner for the system considered is
	given by the block-triangular matrix $\mathcal{P}_1$ given in
	\eqref{optimal_precond}. Again, we are interested in finding approximations
	of $\Phi$ and $S$. For the problems considered here, the matrix $\Phi$ is a
	mass matrix. An efficient way of approximating it is given by a fixed number
	of steps of the Chebyshev semi-iteration
	\cite{GolubVargaI,GolubVargaII,Wathen_Rees}, preconditioned with a Jacobi
	splitting.

We now introduce the matching strategy for approximating the inner Schur complement
	$S = \frac{1}{\beta} \Phi + \Psi_2 \Phi^{-1} \Psi_1^\top$. We seek
	an easily invertible approximation that ``captures'' the two s.p.d.~terms
	in $S$, namely, $\frac{1}{\beta} \Phi$ and $\Psi_2 \Phi^{-1} \Psi_1^\top$.
	The approximation we seek is of the form
	\begin{equation}\label{matching_strategy}
		\widetilde{S} := (\Psi_2 + \Lambda) \Phi^{-1} (\Psi_1 +
			\Lambda)^\top \approx S,
	\end{equation}
	with the matrix $\Lambda$ such that $\Lambda \Phi^{-1} \Lambda^\top =
	\frac{1}{\beta} \Phi$. It is easy to check that, under our assumptions,
	this relation holds if we choose
	\begin{displaymath}
		\Lambda = \dfrac{1}{\sqrt{\beta}} \Phi.
	\end{displaymath}

When the incompressibility constraints are satisfied in the strong sense,
	the approximation given in \eqref{matching_strategy} for the Picard
	iteration is optimal. In fact, in this case, the system is symmetric,
	and each non-linear iteration can be considered as the optimality
	conditions of an incompressible convection--diffusion control problem.
	From here, it can be proved that $\lambda(\widetilde{S}^{-1} S) \subseteq
	[\frac{1}{2}, 1]$; see, for instance, \cite{Pearson_Wathen_2013}.

Remarkably, the approximation given in \eqref{matching_strategy} can also be
	employed when the incompressibility constraints are satisfied only
	inexactly. In fact, one can observe the robustness of this approach when
	employing it as an inner solver for the linear system arising in a Picard
	iteration for the Navier--Stokes control problem;
	see \cite{Leveque_Pearson_2022}. We note that in this case the system is
	non-symmetric, thus the derivation of bounds for the spectrum of the
	preconditioned matrix is non-trivial. Similarly, one can apply the
	matching strategy when solving the linear systems arising from an
	inexact Newton iteration. However, one cannot prove the optimality of
	the strategy. In fact, the upper bound of the spectrum of the preconditioned
	matrix holds if and only if the ``mixed term'' $\Psi_2 + \Psi_1^\top$ is
	s.p.d.; see, for example, \cite{Pearson_Wathen_2013}. Nonetheless, we expect
	this approach to result in a robust preconditioner.

\subsection{Block Pressure Convection--Diffusion preconditioner}\label{sec_3_2} 
We present here the block pressure convection--diffusion preconditioner for the outer Schur complement. This
	heuristic approach has been derived in \cite{Leveque_Pearson_2022} for
	the Navier--Stokes control problem. We would like to mention that
	the authors in \cite{Danieli_Southworth_Wathen} independently derived
	a similar approach for a parallel-in-time solver for the all-at-once
	linear system arising from a Picard iteration applied to the forward
	instationary incompressible Navier--Stokes equations. We follow the
	derivation in \cite{Leveque_Pearson_2022}.

The block pressure convection--diffusion preconditioner is a generalization of
	the pressure convection--diffusion preconditioner derived in
	\cite{Kay_Loghin_Wathen}. Suppose we have a differential operator
	$\mathcal{D}$ of the form
	\begin{displaymath}
		\mathcal{D} = \left[
			\begin{array}{ccc}
				\mathcal{D}^{1,1} & \ldots & \mathcal{D}^{1,n}\\
				\vdots & \ddots & \vdots\\
				\mathcal{D}^{n,1} & \ldots & \mathcal{D}^{n,n}
			\end{array}
		\right],
	\end{displaymath}
	for some $n \in \mathbb{N}$. Here, each $\mathcal{D}^{i,j}$ is a
	differential operator on the velocity space, for $i, j = 1, 2, \ldots, n$.
	Suppose that for each $\mathcal{D}^{i,j}$ the corresponding differential
	operator on the pressure space $\mathcal{D}_p^{i,j}$ is well defined, for
	$i, j = 1, 2, \ldots, n$. Suppose also that the ``commutator" operator
	\begin{displaymath}
		\mathcal{E}_{n} = \mathcal{D} \nabla_{n} - \nabla_{n} \mathcal{D}_p,
	\end{displaymath}
	is small in some sense, where $\nabla_{n} = I_{n} \otimes \nabla$. After
	discretizing with stable finite elements we obtain
	\begin{equation}\label{discretized_general_commutator}
		\left(\mathcal{M}^{-1} \mathbf{D}\right) \mathcal{M}^{-1}
			\vec{B}^{\: \top} - \mathcal{M}^{-1} \vec{B}^{\: \top}
			\left( \mathcal{M}_p^{-1} D_p \right) \approx 0,
	\end{equation}
	where $\mathcal{M}= I_{n} \otimes \mathbf{M}$, $\mathcal{M}_p= I_{n} \otimes
	M_p$, $\vec{B}= I_{n} \otimes B$, and
	\begin{displaymath}
		\mathbf{D} = \left[
			\begin{array}{ccc}
				\mathbf{D}^{1,1} & \ldots & \mathbf{D}^{1,n}\\
				\vdots & \ddots & \vdots\\
				\mathbf{D}^{n,1} & \ldots & \mathbf{D}^{n,n}
			\end{array}
		\right], \qquad
		D_p = \left[
			\begin{array}{ccc}
				D_p^{1,1} & \ldots & D_p^{1,n}\\
				\vdots & \ddots & \vdots\\
				D_p^{n,1} & \ldots & D_p^{n,n}
			\end{array}
		\right].
	\end{displaymath}
	Here, $\mathbf{M}^{-1} \mathbf{D}^{i,j}$ and $M_p^{-1} D^{i,j}_p$ are the
	corresponding discretizations of $\mathcal{D}^{i,j}$ and
	$\mathcal{D}^{i,j}_p$, respectively.
	Assuming invertibility of $\mathbf{D}$ and of $D_p$, by pre-multiplying
	\eqref{discretized_general_commutator} by $\vec{B} \mathbf{D}^{-1}
	\mathcal{M}$ and post-multiplying by $D_p^{-1} \mathcal{M}_p$ we obtain
	\begin{displaymath}
		\vec{B} \mathcal{M}^{-1} \vec{B}^{\: \top} D_p^{-1} \mathcal{M}_p
			\approx \vec{B} \mathbf{D}^{-1} \vec{B}^{\: \top}.
	\end{displaymath}
	We note that $\vec{B} \mathcal{M}^{-1} \vec{B}^{\: \top} = I_{n} \otimes
	(B \mathbf{M}^{-1} B^\top)$. Further, it can be proved that $K_p \approx B
	\mathbf{M}^{-1} B^\top$ for enclosed flow, see
	\cite[pp.\,176--177]{Elman_Silvester_Wathen}. From here, we
	can derive the following approximation:
	\begin{equation}\label{commutator_approx}
		\mathcal{K}_p D_p^{-1} \mathcal{M}_p \approx \vec{B}
			\mathbf{D}^{-1} \vec{B}^{\: \top},
	\end{equation}
	with $\mathcal{K}_p = I_{n} \otimes K_p$. For the problems we consider
	in this work, we have $n=2$, and the matrix $\mathbf{D}$ is the
	$(1,1)$-block of the discretized optimality conditions derived in section
	\ref{sec_2}.

The block pressure convection--diffusion preconditioner has been proved to be
	an efficient and robust approximation of the Schur complement of linear
	systems arising from incompressible flow problems. In fact, despite the
	complex structures of the problems considered in
	\cite{Danieli_Southworth_Wathen, Leveque_Bergamaschi_Martinez_Pearson,
	Leveque_Pearson_2022}, the numerical results obtained showed only
	a mild dependence of the preconditioner with respect to the parameters
	involved. However, this approximation can be applied only to problems
	involving the incompressible Stokes equations or to the Picard approximation
	of problems involving the incompressible Navier--Stokes equations. In
	fact, it is not possible to define a differential operator on the pressure
	space corresponding to the differential operator represented by the matrix
	in \eqref{Newton_derivative_matrix} containing second-order information
	of the convection term. To be specific, the term $\nabla \vec{v}^{\: (k)}$
	is a tensor, while the finite element basis functions on the pressure space
	are scalar functions, thus, if one replaces the finite element basis
	functions on the velocity space in \eqref{Newton_derivative_matrix} with
	the finite element basis functions on the pressure space, one would
	obtain an operator that is not well defined over the pressure
	space.
	For this reason, one should not
	expect this approach to result in a robust solver when applied to a
	Newton iteration. In the following section, we introduce an approach
	that obviates this issue.

\subsection{Augmented Lagrangian preconditioner}\label{sec_3_3}
In this section we develop an augmented Lagrangian preconditioner for the outer Schur complement, as an alternative to the block pressure convection--diffusion approach of the previous section.
As mentioned above, the most delicate task when deriving a robust
	preconditioner of the form given in \eqref{approx_preconditioners} is
	finding a good approximation $\widetilde{S}$ of the Schur complement $S$.
	A powerful approach that circumvents this issue is given by the augmented
	Lagrangian formulation \cite{Fortin_Glowinski}. An additional benefit 
	is that, when applied to incompressible flow problems, this
	approach can be regarded as a form of grad-div stabilization, in which
	an extra term is added to the momentum equation in order to better enforce
	the incompressibility of the fluid, see
	\cite{Olshanskii_Lube_Heister_Lowe, Olshanskii_Reusken}.
	The extra term allows for improved stability and accuracy of some
	discretizations, and leads to the construction of robust preconditioning
	techniques

Augmented Lagrangian-based preconditioning has been shown to be very effective
	for solving linearizations of incompressible fluid flow problems. One such 
	technique was developed in \cite{Benzi_Olshanskii_2006} for the
	solution of the Picard iteration applied to the (forward) incompressible
	Navier--Stokes equations in 2D, and has been successfully extended to the
	Newton iteration of the 3D incompressible Navier--Stokes equations
	\cite{Farrell_Mitchell_Wechsung}, to problems with indefinite $(1,1)$-block
	\cite{Benzi_Liu,Olshanskii_Benzi}, as well as to other problems of
	saddle-point type. We also mention the recent work \cite{Lohmann_Turek},
	where the authors derived an augmented Lagrangian preconditioner for the
	all-at-once linear system arising upon discretization of the instationary
	incompressible Navier--Stokes equations, when employing a Picard
	linearization.

We begin presenting the idea of the augmented Lagrangian preconditioner
	for general saddle-point systems, giving as an example the linearization of the
	stationary incompressible Navier--Stokes equations. Then, we specialize the
	framework to the case of linearizations of the stationary Navier--Stokes
	control problem considered in this work.

Given the linear system in \eqref{generalized_saddle_point}, with $\Theta = 0$
	and $\Psi_1 = \Psi_2=\Psi$, the idea of the augmented Lagrangian
	preconditioner is to consider an equivalent system of the form
	\begin{equation}\label{augmented_system}
		\underbrace{\left[
			\begin{array}{cc}
				\Phi + \gamma \Psi^\top \mathcal{W}^{-1} \Psi
					& \Psi^\top\\
				\Psi & 0
			\end{array}
		\right]}_{\mathcal{A}_{\gamma}}
		\left[
			\begin{array}{c}
				\bm{y}_1\\
				\bm{y}_2
			\end{array}
		\right]
		=
		\left[
			\begin{array}{c}
				\bm{\hat{b}}_1\\
				\bm{b}_2
				\end{array}
		\right],
	\end{equation}
	for suitable $\gamma>0$ and $\mathcal{W}$, where we set $\bm{\hat{b}}_1 =
	\bm{b}_1 + \gamma \Psi_1^\top \mathcal{W}^{-1} \bm{b}_2$. For
	example, in the setting of the incompressible
	Navier--Stokes equations, for which we have $\Psi
	= B$, the matrix $\mathcal{W}$ is chosen to be the pressure
	mass matrix $M_p$ or its diagonal. On the other hand, the choice of
	$\gamma$ is more delicate, as we describe below.

In order to solve the system in \eqref{augmented_system}, we employ the 
	block-triangular preconditioner $\mathcal{P}_2$ given in
	\eqref{optimal_precond}. Given the structure of the $(1,1)$-block
	of $\mathcal{A}_\gamma$ in \eqref{augmented_system}, with
	a small change in notation\footnote{Note that
	the preconditioner $\mathcal{P}_2$ specializes to $\mathcal{P}_\gamma$ for
	the augmented Lagrangian approach.}, we consider the ideal preconditioner
	\begin{equation}\label{augmented_precond}
		\mathcal{P}_{\gamma} =
			\left[
				\begin{array}{cc}
					\Phi + \gamma \Psi^\top \mathcal{W}^{-1} \Psi
						& \Psi^\top\\
					0 & - S_{\gamma}
				\end{array}
			\right],
	\end{equation}
	where $S_{\gamma}=\Psi (\Phi + \gamma \Psi^\top \mathcal{W}^{-1}
	\Psi)^{-1} \Psi^\top$. Assuming that the matrix $\Psi \Phi^{-1}
	\Psi^\top$ is invertible and employing the Sherman--Morrison--Woodbury
	formula, one can show that
	\begin{equation}\label{S_gamma_inverse}
		S_{\gamma}^{-1} = (\Psi \Phi^{-1} \Psi^\top)^{-1} +
			\gamma \mathcal{W}^{-1}.
	\end{equation}
	In fact, in this case the Sherman--Morrison--Woodbury formula writes as
	\begin{displaymath}
		(\Phi + \gamma \Psi^\top \mathcal{W}^{-1} \Psi)^{-1} =
			\Phi^{-1} - \Phi^{-1} \Psi^\top \left( \gamma^{-1} \mathcal{W}
			+ \Psi \Phi^{-1} \Psi^\top \right)^{-1} \Psi \Phi^{-1},
	\end{displaymath}
	which, together with the invertibility of $\Psi \Phi^{-1} \Psi^\top$,
	implies that
	\begin{align}
		S_{\gamma} & = \Psi \Phi^{-1} \Psi^\top - \Psi \Phi^{-1} \Psi^\top
			\left( \gamma^{-1} \mathcal{W} +
			\Psi \Phi^{-1} \Psi^\top \right)^{-1} \Psi \Phi^{-1} \Psi^\top
			\nonumber \\
			& = \Psi \Phi^{-1} \Psi^\top \left( \gamma^{-1} \mathcal{W} +
			\Psi \Phi^{-1} \Psi^\top \right)^{-1} \left( \gamma^{-1}
			\mathcal{W} + \Psi \Phi^{-1} \Psi^\top - \Psi \Phi^{-1} \Psi^\top
			\right) \nonumber \\
			& = \left[ (\Psi \Phi^{-1} \Psi^\top)^{-1} + \gamma \mathcal{W}^{-1} 
			\right]^{-1}. \nonumber
	\end{align}
	From the derivation above, we understand that an approximate
	inverse of $S_{\gamma}$ can be obtained given a ``rough'' approximation of
	$(\Psi \Phi^{-1} \Psi^\top)^{-1}$. However, above all the choice
	of the parameter $\gamma$ is important in order to obtain a fast solver.
	For example, in the case of 
	saddle-point systems  arising from the linearization and discretization
	of the stationary Navier--Stokes equations, the term $\Psi \Phi^{-1}
	\Psi^\top$ is replaced by a multiple of the mass matrix $M_p$ (or its
	diagonal). In particular, the following approximation is employed:
	\begin{displaymath}
		\widetilde{S}_\gamma ^{-1} = \nu M_p^{-1} + \gamma M_p^{-1},
			\quad \mathrm{or} \quad \widetilde{S}_\gamma =
			\dfrac{1}{\nu + \gamma} M_p.
	\end{displaymath}
	In \cite{Benzi_Olshanskii_2011}, the authors proved that the choice of
	$\gamma \sim \nu^{-1}$ results in an $h$- and
	$\nu$-independent preconditioner, but numerical
	experiments show that even smaller values of $\gamma$ result in a
	robust solver, see \cite{Benzi_Olshanskii_2006}.

We now move to the derivation of the preconditioner for the Navier--Stokes
	control problem. In order to simplify the notations, we consider the
	linear system \eqref{Opt_then_Discr_discrete_Inexact_Newton_system} arising
	from the discretization of the inexact Newton iteration by employing an
	optimize-then-discretize approach, noting that the derivation easily
	generalizes to the discretize-then-optimize one.

As mentioned above, the first step consists of considering an equivalent
	system of the form
	\begin{displaymath}
		\underbrace{\left[
			\begin{array}{cc}
				{\mathcal F}^{(k)}_{\mathrm{OD}} + \gamma
					{\mathcal B}^\top \mathcal{W}^{-1} {\mathcal B} &
					{\mathcal B}^\top \\
				{\mathcal B} & 0
			\end{array}
		\right]}_{\mathcal{A}_{\gamma}^{(k)}}
		\left[
			\begin{array}{c}
				\bm{\delta v}^{(k)}\\
				\bm{\delta \zeta}^{(k)}\\
				\bm{\delta \mu}^{(k)}\\
				\bm{\delta p}^{(k)}
			\end{array}
		\right]=
		\left[
			\begin{array}{c}
				\bm{\hat{R}}_1^{(k)}\\
				\bm{\hat{R}}_2^{(k)}\\
				\bm{r}_1^{(k)}\\
				\bm{r}_2^{(k)}
			\end{array}
		\right],
	\end{displaymath}
	for $\gamma>0$ and a suitable matrix $\mathcal{W}$. Here, the
	vectors $\bm{\hat{R}}_1^{(k)}$ and $\bm{\hat{R}}_2^{(k)}$ are given by
	\begin{displaymath}
		\left[
			\begin{array}{cc}
				\bm{\hat{R}}_1^{(k)}\\
				\bm{\hat{R}}_2^{(k)}
			\end{array}
		\right]
		=
		\left[
			\begin{array}{cc}
				\bm{R}_1^{(k)}\\
				\bm{R}_2^{(k)}
			\end{array}
		\right] +
		\gamma \Psi^\top \mathcal{W}^{-1}
		\left[
			\begin{array}{cc}
				\bm{r}_1^{(k)}\\
				\bm{r}_2^{(k)}
			\end{array}
		\right].
	\end{displaymath}
	We employ as matrix $\mathcal{W}$ the following matrix:
	\begin{equation}\label{mathcal_W}
		\mathcal{W} =
		\left[
			\begin{array}{cc}
				0 & W\\
				W & 0
			\end{array}
		\right],
	\end{equation}
	with $W$ being the pressure mass matrix $M_p$ or its diagonal. With this
	choice of $\mathcal{W}$, the augmented $(1,1)$-block becomes
	\begin{displaymath}
	{\mathcal F}^{(k)}_{\mathrm{OD}} + \gamma {\mathcal B}^\top
		\mathcal{W}^{-1} {\mathcal B} = \left[
			\begin{array}{cc}
				\mathbf{M} &
					\mathbf{D}^{(k)}_{\mathrm{adj}} + (\mathbf{H}^{(k)})^\top
						+ \gamma B^\top W^{-1} B\\
				\mathbf{D}^{(k)} + \mathbf{H}^{(k)}
					+ \gamma B^\top W^{-1} B & -\mathbf{M}_\beta
			\end{array}
		\right].
	\end{displaymath}
	We note that this choice of $\mathcal{W}$ can be regarded
	as a grad-div stabilization of both the state and adjoint equations. This
	stabilization can be derived in a formal way starting from the control
	problem \eqref{Navier_Stokes_control_cost_functional}--\eqref{Navier_Stokes_control_constraints}, adding the grad-div stabilization
	term in the constraints \eqref{Navier_Stokes_control_constraints}, and
	finally deriving the optimality conditions for this problem, as done in
	\cite{Cibik}.

We end this section by finding a suitable approximation of the Schur complement
	$S_{\gamma}$ in \eqref{augmented_precond}. By applying the
	Sherman--Morrison--Woodbury formula, one obtains that the inverse of
	$S_{\gamma}$ can be written as in \eqref{S_gamma_inverse}. In our case, the
	matrix ${\mathcal F}^{(k)}_{\mathrm{OD}}$ depends on the regularization
	parameter $\beta$. In numerical experiments, we observed that the term
	$({\mathcal B} ({\mathcal F}^{(k)}_{\mathrm{OD}})^{-1}
	{\mathcal B}^\top)^{-1}$ in \eqref{S_gamma_inverse} is not negligible
	for small $\beta$. Therefore, in order to obtain a robust solver, we need
	to find a suitable approximation of it. Thus, in order to improve the
	robustness of the solver, we consider the following approximation
	of $S_{\gamma}^{-1}$:
	\begin{equation}\label{widetilde_S_gamma_inverse}
		\widetilde{S}_{\gamma}^{-1} = 
			\left[
				\begin{array}{cc}
					K_p^{-1} & \gamma W^{-1}\\
					\gamma W^{-1} & -\dfrac{1}{\beta} K_p^{-1}
				\end{array}
			\right].
	\end{equation}
	We made this choice as, for small $\beta$, the matrix
	${\mathcal F}^{(k)}_{\mathrm{OD}}$ is well approximated by a multiple of
	a block-diagonal mass matrix. Specifically, for $\beta$ small,
	we have
	\begin{displaymath}
	{\mathcal F}^{(k)}_{\mathrm{OD}} \approx
		\left[
			\begin{array}{cc}
				\mathbf{M} & 0\\
				0 & - \dfrac{1}{\beta} \mathbf{M}
			\end{array}
		\right].
	\end{displaymath}
	Then, using the fact that $K_p \approx B \mathbf{M}^{-1} B^\top$, we
	obtain the proposed approximation. Finally, the
	observation above gives us a heuristic choice of the parameter $\gamma$. In
	fact, the matrix ${\mathcal F}^{(k)}_{\mathrm{OD}} + \gamma \Psi^\top
	\mathcal{W}^{-1} \Psi$ is similar to the following matrix:
	\begin{displaymath}
		 \underbrace{
			\left[
				\begin{array}{cc}
					\mathbf{M} &
					\sqrt{\beta}(\mathbf{D}^{(k)}_{\mathrm{adj}} +
						(\mathbf{H}^{(k)})^\top) \\
					\sqrt{\beta}(\mathbf{D}^{(k)} + \mathbf{H}^{(k)}) &
					-\mathbf{M}
				\end{array}
			\right]
		}_{\widehat{\mathcal F}^{(k)}_{\mathrm{OD}}} + \,
		\bar{\gamma}
		\left[
			\begin{array}{cc}
				0 & B^\top W^{-1} B \\
				B^\top W^{-1} B & 0
			\end{array}
		\right],
	\end{displaymath}
	with $\bar{\gamma}= \sqrt{\beta}\gamma$. Following a similar argument
	as in \cite{Benzi_Olshanskii_2011}, one can expect that for large
	values of $\bar{\gamma}$ the augmented Lagrangian preconditioner is able to
	cluster the majority of the eigenvalues of the preconditioned matrix around
	$1$. This is because the term
	$\bar{\gamma} {\mathcal B}^\top \mathcal{W}^{-1} {\mathcal B}$
	is becoming dominant with respect to the matrix
	$\widehat{\mathcal F}^{(k)}_{\mathrm{OD}}$, thus the way one approximates
	the inverse of the matrix $\Psi \Phi^{-1} \Psi^\top$ in
	\eqref{S_gamma_inverse} becomes less critical for large
	$\bar{\gamma}$. Finally, the previous observation leads to
	the heuristic choice of $\gamma > \frac{1}{\sqrt{\beta}}$. 
	A more formal eigenvalue analysis is beyond the scope of
	this work, and will be the topic of future research.

Finally, as we mentioned above, practical preconditioners usually
	require inexact solves for the main blocks. For this reason, the preconditioner
	we employ in the numerical tests is an approximation
	$\widetilde{\mathcal{P}}_\gamma$ of the ideal preconditioner
	$\mathcal{P}_\gamma$ given in \eqref{augmented_precond}. Specifically, we
	approximate the $(1,1)$-block ${\mathcal F}^{(k)}_{\mathrm{OD}} +
	\gamma {\mathcal B}^\top \mathcal{W}^{-1} {\mathcal B}$ with a fixed
	number of iterations of the FGMRES solver preconditioned by a
	suitable preconditioner, and employ the matrix $\widetilde{S}_{\gamma}$ in
	\eqref{widetilde_S_gamma_inverse} as an approximation of the Schur
	complement $S_\gamma$.

\section{Numerical results}\label{sec_4}
In this section, we present numerical evidence of the robustness of the
	augmented Lagrangian approach.

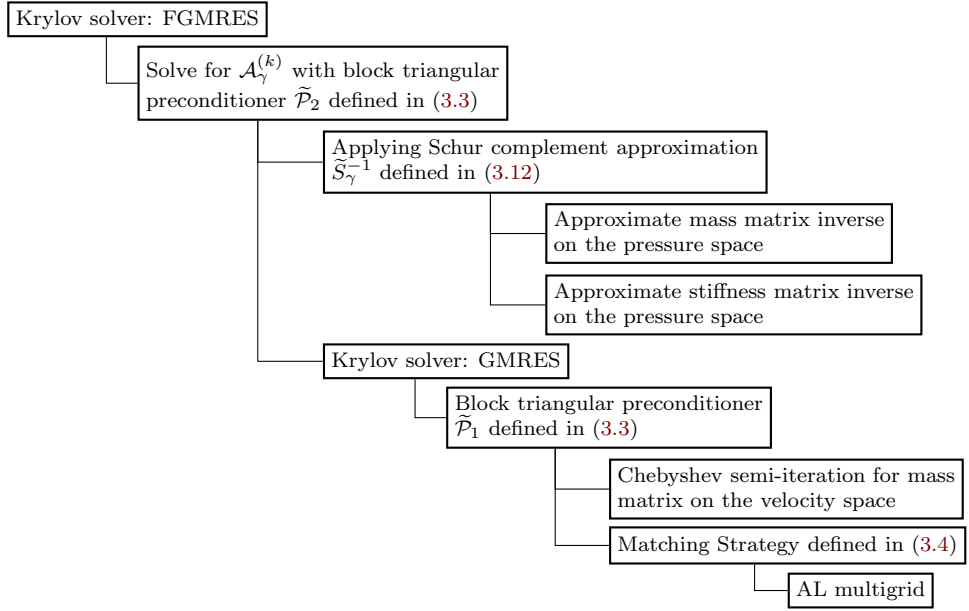
\begin{figure}[tbhp]
	\footnotesize
	\centering
	\begin{tikzpicture}[%
		every node/.style={draw=black, thick, anchor=west},
		grow via three points={one child at (-0.0,-0.7) and two children at (0.0,-0.7) and (0.0,-1.4)},
		edge from parent path={(\tikzparentnode.210) |- (\tikzchildnode.west)}]
		\node {Krylov solver: FGMRES}
		child {node at (0.0,-0.15) [align=left] {Solve for
			$\mathcal{A}_\gamma^{(k)}$ with block triangular\\
			preconditioner $\widetilde{\mathcal{P}}_\gamma$ defined in
			Section \ref{sec_3_3}}
			child {node at (0.0,-0.35) [align=left] {Applying Schur complement
				approximation\\ $\widetilde{S}_\gamma^{-1}$ defined in
				\eqref{widetilde_S_gamma_inverse}}
				child {node at (0.0,-0.25) [align=left] {Approximate mass matrix inverse\\
					on the pressure space}}
				child {node at (0.0,-0.5) [align=left]{Approximate stiffness matrix inverse\\
					on the pressure space}}
			}
			child {node at (0.0,-2.3) {Krylov solver: GMRES}
				child { node at (0.0,-0.05) [align=left] {Block triangular
					preconditioner\\
					$\widetilde{\mathcal{P}}_1$ defined in
					\eqref{approx_preconditioners}}
				child {node at (0.0,-0.25) [align=left] {Chebyshev
					semi-iteration for mass\\
					matrix on the velocity space}}
				child {node at (0.0,-0.3) {Matching Strategy defined in
					\eqref{matching_strategy}}
					child {node at (0.0,0.1) {AL multigrid}}
					child[missing]{}
					child[missing]{}
				}
				child[missing]{}
				}
			}
			child[missing]{}
			child[missing]{}
			child[missing]{}
		};
	\end{tikzpicture}
	\caption{Solver diagram for the augmented Lagrangian preconditioner applied
	to stationary Navier--Stokes control problems.}
	\label{fig:AL_NavierStokes_Control}
	\end{figure}

In all our tests, we employ as a preconditioner the block-triangular matrix
	$\widetilde{\mathcal{P}}_2$ given in \eqref{approx_preconditioners}.
	Note that, as we mentioned above, in the augmented Lagrangian
	framework the preconditioner $\widetilde{\mathcal{P}}_2$ specializes to
	$\widetilde{\mathcal{P}}_\gamma$ as described at the end of
	Section \ref{sec_3_3}. The inverse of the $(1,1)$-block that couples the two
	velocities is approximated by 5 iterations of
	preconditioned GMRES \cite{Saad_Schultz}, employing the block-triangular
	matrix $\widetilde{\mathcal{P}}_1$ defined in \eqref{approx_preconditioners}
	as a preconditioner. In the inner solver, we approximate the inverse of the
	$(1,1)$-block with 20 steps of Chebyshev semi-iteration, while the Schur
	complement is approximated by the matrix $\widetilde{S}$ defined in
	\eqref{matching_strategy}. Since we need an inner Krylov solver and our
	systems are non-symmetric, we employ flexible GMRES \cite{Saad} as outer
	solver, with restart every 10 iterations, up to a tolerance $10^{-6}$ on
	the relative residual. All CPU times below are reported in seconds.
	We summarise the preconditioner for the augmented Lagrangian method
	in Figure \ref{fig:AL_NavierStokes_Control}. Following our observation
	at the end of Section \ref{sec_3_3}, for a given $\beta$ we set
	$\gamma=10/\sqrt{\beta}$.

In all our tests, we report the value of the discrete cost functional
	$J_h(\bm{v}^*, \bm{u}^*)$ at the optimal solution $(\bm{v}^*, \bm{u}^*)$ found
	at the end of the non-linear process. The discrete cost functional is defined as
	follows:
	\begin{displaymath}
	J^*_h:= J_h(\bm{v}^*, \bm{u}^*) = \dfrac{1}{2}(\bm{v}^* - \bm{v}_d)^\top
		\mathbf{M} (\bm{v}^* - \bm{v}_d) + \dfrac{\beta}{2}(\bm{u}^*)^\top
		\mathbf{M} \bm{u}^*.
	\end{displaymath}

\subsection{Lid-Driven Cavity}\label{sec_4_1}
For the first test we consider here, we set $d=2$ (that is, $\mathbf{x}=(x_1,x_1)$)
	and $\Omega=(-1,1)^2$. The Navier--Stokes control problem is given
	by \eqref{Navier_Stokes_control_cost_functional}--\eqref{Navier_Stokes_control_constraints}, with $\vec{f}=\vec{v}_d=\vec{0}$,
	and
	\begin{displaymath}
		\vec{g} = \left\{
			\begin{array}{ll}
				\left[1,0\right]^\top & \mathrm{on} \: \partial \Omega_1
					:= \left(-1,1 \right)\times \left\{1\right\},\\
				\left[0,0\right]^\top & \mathrm{on} \:\partial \Omega
					\setminus \partial\Omega_1.
			\end{array}
		\right.
	\end{displaymath}

We employ both the optimize-then-discretize and the discretize-then-optimize
	strategy, showing that the augmented Lagrangian approach can obtain fast
	convergence for solving the Navier--Stokes control problem, exhibiting
	a robust behaviour with respect to problem parameters. In the first setting, the
	optimize-then-discretize one, we compare the augmented Lagrangian preconditioner
	with the block pressure convection--diffusion preconditioner; this is done
	by applying the inverses of the velocity blocks in \eqref{matching_strategy}
	exactly, since for an inexact solver the augmented velocity blocks require
	a specific geometric multigrid not easily available in MATLAB. Then, in
	the discretize-then-optimize setting, we give as a proof of concept that
	the whole augmented Lagrangian preconditioner can be run in an inexact
	framework, by applying a fixed number of cycles of a specific geometric
	multigrid on the augmented velocity blocks built into Firedrake. The results
	obtained from the second test allow us to observe the strength of the
	augmented Lagrangian preconditioner, that allows one to obtain a prescribed
	reduction on the relative linear residual in very few iterations and, more
	importantly, independently of the problem parameters even if the
	preconditioner is applied inexactly.

We would like to note here that for this problem the desired state does not solve
	the stationary incompressible Navier--Stokes equations problem with the imposed
	boundary conditions. This is due to a mismatch of the boundary conditions on
	the ``lid'' of the square domain. For this reason, we cannot expect the cost
	functional to be driven to $0$ as the mesh is refined.

\subsubsection{Optimize-then-Discretize}\label{sec_4_1_1}
We first compare the augmented Lagrangian strategy with the block
	pressure convection--diffusion preconditioner, in an
	optimize-then-discretize framework. Specifically, we are interested in
	comparing the two Schur complement approximations. The tests presented here
	are run on MATLAB R2018b, using a 1.70GHz Intel quad-core i5 processor and
	8 GB RAM on an Ubuntu 18.04.1 LTS operating system.

The outer solver is based on the flexible GMRES routine in the TT-Toolbox
	\cite{TT_Toolbox}, while we employ the GMRES routine implemented in MATLAB
	as inner solver. The inverse of the stiffness matrix on the pressure space
	is approximated by 2 V-cycles (with 2 symmetric Gauss--Seidel iterations
	for pre-/post-smoothing) of the \texttt{HSL\_MI20} solver \cite{HSL_MI20}.
	We apply 20 steps of Chebyshev semi-iteration as an approximate inverse of
	the mass matrix on the pressure space. Since we are interested in
	understanding the behaviour of the Schur complement approximations,
	we solve for the approximation given in \eqref{matching_strategy}
	exactly, that is, the blocks $\Psi_2 + \Lambda$ and
	$(\Psi_1 + \Lambda)^\top$ are solved by MATLAB backslash. We employ the
	diagonal of the mass matrix on the pressure space as matrix $W$ in
	\eqref{mathcal_W}. Finally, when applying the block pressure
	convection--diffusion preconditioner, we approximate the Schur complement
	with the matrix defined in \eqref{commutator_approx} (the differential
	operator on the pressure space does not contain the Newton matrices),
	while we employ the approximation given in
	\eqref{widetilde_S_gamma_inverse} for the augmented Lagrangian approach.

Regarding the non-linear iteration, we solve each problem up
	to a reduction of $10^{-5}$ on the non-linear relative residual. We allow
	10 iterations for inexact Newton. The first non-linear iteration is employed
	for solving the corresponding Stokes control problem.

We use inf--sup stable Taylor--Hood $[{Q}_2]^2$--${Q}_1$ finite
	elements in the spatial dimensions, with level of refinement $l$
	representing a (spatial) uniform grid of mesh-size $h=2^{1-l}$ for
	${Q}_1$ basis functions, and $h=2^{-l}$ for ${Q}_2$ elements,
	in each dimension. For this test case, we employ the
	LPS approach for stabilizing the linearization of the Navier--Stokes
	control problem at each Newton iteration.

\begin{table}[!ht]
\caption{Average GMRES iterations with the block pressure convection--diffusion preconditioner, for $\nu=\frac{1}{100}$, $\frac{1}{250}$, and $\frac{1}{500}$, and a range of $l$, $\beta$.}\label{Inexact_Newton_Opt_then_Discr_BPCD}
\begin{scriptsize}
\begin{center}
\renewcommand{\arraystretch}{1.2}
{\begin{tabular}{|c||c|c|c|c|c||c|c|c|c|c||c|c|c|c|c|}
\hline
\multicolumn{1}{|c||}{} & \multicolumn{5}{c||}{$\nu=\frac{1}{100}$} & \multicolumn{5}{c||}{$\nu=\frac{1}{250}$} & \multicolumn{5}{c|}{$\nu=\frac{1}{500}$}\\
\cline{2-16}
\multicolumn{1}{|c||}{} & \multicolumn{5}{c||}{$\beta$} & \multicolumn{5}{c||}{$\beta$} & \multicolumn{5}{c|}{$\beta$}\\
\cline{2-16}
$l$ & \!\!\!\! $10^{-1}$ \!\!\!\! & \!\!\!\! $10^{-2}$ \!\!\!\! & \!\!\!\! $10^{-3}$ \!\!\!\! & \!\!\!\! $10^{-4}$ \!\!\!\! & \!\!\!\! $10^{-5}$ \!\!\!\! & \!\!\!\! $10^{-1}$ \!\!\!\! & \!\!\!\! $10^{-2}$ \!\!\!\! & \!\!\!\! $10^{-3}$ \!\!\!\! & \!\!\!\! $10^{-4}$ \!\!\!\! & \!\!\!\! $10^{-5}$ \!\!\!\! & \!\!\!\! $10^{-1}$ \!\!\!\! & \!\!\!\! $10^{-2}$ \!\!\!\! & \!\!\!\! $10^{-3}$ \!\!\!\! & \!\!\!\! $10^{-4}$ \!\!\!\! & \!\!\!\! $10^{-5}$ \!\!\!\! \\
\hline
\hline
\!\!\!\! $3$ \!\!\!\! & \!\!\!\! 30 \!\!\!\! & \!\!\!\! 16 \!\!\!\! & \!\!\!\! 11 \!\!\!\! & \!\!\!\! 9 \!\!\!\! & \!\!\!\! 7 \!\!\!\! & \!\!\!\! 33 \!\!\!\! & \!\!\!\! 15 \!\!\!\! & \!\!\!\! 10 \!\!\!\! & \!\!\!\! 9 \!\!\!\! & \!\!\!\! 7 \!\!\!\! & \!\!\!\! 31 \!\!\!\! & \!\!\!\! 15 \!\!\!\! & \!\!\!\! 10 \!\!\!\! & \!\!\!\! 9 \!\!\!\! & \!\!\!\! 7 \!\!\!\! \\
\hline
\!\!\!\! $4$ \!\!\!\! & \!\!\!\! 36 \!\!\!\! & \!\!\!\! 21 \!\!\!\! & \!\!\!\! 13 \!\!\!\! & \!\!\!\! 10 \!\!\!\! & \!\!\!\! 8 \!\!\!\! & \!\!\!\! 64 \!\!\!\! & \!\!\!\! 22 \!\!\!\! & \!\!\!\! 12 \!\!\!\! & \!\!\!\! 10 \!\!\!\! & \!\!\!\! 8 \!\!\!\! & \!\!\!\! 75 \!\!\!\! & \!\!\!\! 22 \!\!\!\! & \!\!\!\! 12 \!\!\!\! & \!\!\!\! 10 \!\!\!\! & \!\!\!\! 8 \!\!\!\! \\
\hline
\!\!\!\! $5$ \!\!\!\! & \!\!\!\! 39 \!\!\!\! & \!\!\!\! 21 \!\!\!\! & \!\!\!\! 16 \!\!\!\! & \!\!\!\! 12 \!\!\!\! & \!\!\!\! 9 \!\!\!\! & \!\!\!\! 75 \!\!\!\! & \!\!\!\! 31 \!\!\!\! & \!\!\!\! 16 \!\!\!\! & \!\!\!\! 11 \!\!\!\! & \!\!\!\! 9 \!\!\!\! & \!\!\!\! 150 \!\!\!\! & \!\!\!\! 35 \!\!\!\! & \!\!\!\! 16 \!\!\!\! & \!\!\!\! 11 \!\!\!\! & \!\!\!\! 9 \!\!\!\! \\
\hline
\!\!\!\! $6$ \!\!\!\! & \!\!\!\! 41 \!\!\!\! & \!\!\!\! 22 \!\!\!\! & \!\!\!\! 17 \!\!\!\! & \!\!\!\! 14 \!\!\!\! & \!\!\!\! 11 \!\!\!\! & \!\!\!\! 73 \!\!\!\! & \!\!\!\! 29 \!\!\!\! & \!\!\!\! 19 \!\!\!\! & \!\!\!\! 14 \!\!\!\! & \!\!\!\! 10 \!\!\!\! & \!\!\!\! 131 \!\!\!\! & \!\!\!\! 41 \!\!\!\! & \!\!\!\! 22 \!\!\!\! & \!\!\!\! 13 \!\!\!\! & \!\!\!\! 10 \!\!\!\! \\
\hline
\!\!\!\! $7$ \!\!\!\! & \!\!\!\! 38 \!\!\!\! & \!\!\!\! 25 \!\!\!\! & \!\!\!\! 19 \!\!\!\! & \!\!\!\! 16 \!\!\!\! & \!\!\!\! 13 \!\!\!\! & \!\!\!\! 64 \!\!\!\! & \!\!\!\! 31 \!\!\!\! & \!\!\!\! 21 \!\!\!\! & \!\!\!\! 16 \!\!\!\! & \!\!\!\! 13 \!\!\!\! & \!\!\!\! 140$\dagger$\footnotemark \!\!\!\! & \!\!\!\! 24 \!\!\!\! & \!\!\!\! 25 \!\!\!\! & \!\!\!\! 17 \!\!\!\! & \!\!\!\! 13 \!\!\!\! \\
\hline
\end{tabular}}
\end{center}
\end{scriptsize}
\end{table}

In Table \ref{Inexact_Newton_Opt_then_Discr_BPCD}, we report the average number
	of flexible GMRES (FGMRES) iterations (rounded to the nearest integer)
	 required for solving the inexact Newton iteration, when
	employing the block pressure convection--diffusion preconditioner for
	approximating the Schur complement. 
	Further, in Tables
	\ref{Inexact_Newton_Opt_then_Discr_Augmented_Lagrangian_1}--\ref{Inexact_Newton_Opt_then_Discr_Augmented_Lagrangian_3} we report
	the average number of FGMRES iterations required for solving the
	corresponding problem, when employing the augmented Lagrangian
	preconditioner, for our heuristic choice of $\gamma$, together with the value
	of the discrete cost functional $J^*_h$. We report the value of the cost
	functional $J^*_h$ only for the augmented Lagrangian approach, as this
	coincides up to the second significant digit with the one obtained by
	applying the block pressure convection--diffusion preconditioner.
	Finally, in Table
	\ref{Number_Inexact_newton_Opt_then_Discr} we report the total number of
	inexact Newton iterations required for obtaining convergence (when applying
	the augmented Lagrangian preconditioner), together with the total degrees
	of freedom (DoF) of the system solved at each non-linear iteration.
	 Again, all the linear solves involved
	in the application of the two preconditioners are performed exactly by backslash.

\footnotetext{$\dagger$ means that the outer (inexact Newton) iteration did not converge in 10 iterations. The average number of FGMRES iterations is evaluated over the first 5 inexact Newton iterations.}

From Table \ref{Inexact_Newton_Opt_then_Discr_BPCD}, we observe the
	$h$-robustness of the block pressure convection--diffusion preconditioner
	 for small values of $\beta$.  For larger values of $\beta$, 
	however, the preconditioner is strongly dependent on the viscosity $\nu$ and 
	lacks robustness. In
	fact, the linear solver is not able to converge for small values of
	the viscosity and larger $\beta$, resulting in a non-convergent
	inexact Newton method. On the other hand, as observed from Tables
	\ref{Inexact_Newton_Opt_then_Discr_Augmented_Lagrangian_1}--\ref{Inexact_Newton_Opt_then_Discr_Augmented_Lagrangian_3}
	the augmented Lagrangian preconditioner is able to reach convergence in
	an average number of linear iterations that is roughly constant with
	respect to all the parameters involved in the problems, especially for
	sufficiently fine grids. In particular,
	we observe that the augmented Lagrangian preconditioner requires at most
	9 iterations in average for reaching the prescribed tolerance on the
	relative residual. The efficiency of the approach is not lost even
	if one drives the regularization parameter towards $0$, converging in at most
	5 iterations for $\beta=10^{-6}$ and $10^{-8}$. Further, from Tables \ref{Inexact_Newton_Opt_then_Discr_Augmented_Lagrangian_1}--\ref{Inexact_Newton_Opt_then_Discr_Augmented_Lagrangian_3} we observe that the
	discrete cost functional $J^*_h$ is decreasing as the mesh size $h$ goes to $0$.
	Finally, from Table
	\ref{Number_Inexact_newton_Opt_then_Discr} we note that the number of
	inexact Newton iterations required for reaching a prescribed non-linear
	tolerance is roughly constant with respect to the regularization
	parameter and the viscosity, for fine grids, with a small dependence on
	the viscosity for coarse grids.

\begin{table}[!ht]
\caption{Average GMRES iterations and value of the discrete cost functional $J^*_h$ with the augmented Lagrangian preconditioner with  $\gamma=10/\sqrt{\beta}$, for $\nu=\frac{1}{100}$ and a range of $l$, $\beta$.}\label{Inexact_Newton_Opt_then_Discr_Augmented_Lagrangian_1}
\begin{scriptsize}
\begin{center}
\renewcommand{\arraystretch}{1.2}
{\begin{tabular}{|c||c|c||c|c||c|c||c|c||c|c||c|c||c|c|}
\hline
\multicolumn{1}{|c||}{} & \multicolumn{14}{c|}{$\beta$} \\
\cline{2-15}
\multicolumn{1}{|c||}{} & \multicolumn{2}{c||}{$10^{-1}$} & \multicolumn{2}{c||}{$10^{-2}$} & \multicolumn{2}{c||}{$10^{-3}$} & \multicolumn{2}{c||}{$10^{-4}$} & \multicolumn{2}{c||}{$10^{-5}$} & \multicolumn{2}{c||}{$10^{-6}$} & \multicolumn{2}{c|}{$10^{-8}$} \\
\cline{2-15}
$l$ & $\mathrm{it}$ & $J^*_h$ & $\mathrm{it}$ & $J^*_h$ & $\mathrm{it}$ & $J^*_h$ & $\mathrm{it}$ & $J^*_h$ & $\mathrm{it}$ & $J^*_h$ & $\mathrm{it}$ & $J^*_h$ & $\mathrm{it}$ & $J^*_h$\\
\hline
\hline
\!\!\!\! $3$ \!\!\!\! & 6 & 4.2e-02 & 5 & 3.3e-02 & 4 & 3.1e-02 & 3 & 3.1e-02 & 3 & 3.0e-02 & 3 & 3.0e-02 & 2 & 3.0e-02 \\
\hline
\!\!\!\! $4$ \!\!\!\! & 7 & 3.7e-02 & 6 & 2.3e-02 & 5 & 1.7e-02 & 4 & 1.5e-02 & 4 & 1.5e-02 & 3 & 1.5e-02 & 3 & 1.5e-02 \\
\hline
\!\!\!\! $5$ \!\!\!\! & 7 & 3.6e-02 & 5 & 2.1e-02 & 5 & 1.3e-02 & 6 & 8.8e-03 & 4 & 7.7e-03 & 4 & 7.5e-03 & 3 & 7.5e-03 \\
\hline
\!\!\!\! $6$ \!\!\!\! & 6 & 3.6e-02 & 6 & 2.1e-02 & 5 & 1.2e-02 & 8 & 7.0e-03 & 5 & 4.6e-03 & 5 & 3.9e-03 & 3 & 3.7e-03 \\
\hline
\!\!\!\! $7$ \!\!\!\! & 6 & 3.6e-02 & 5 & 2.1e-02 & 5 & 1.2e-02 & 6 & 6.9e-03 & 7 & 3.9e-03 & 5 & 2.5e-03 & 4 & 1.9e-03 \\
\hline
\end{tabular}}
\end{center}
\end{scriptsize}
\end{table}

\begin{table}[!ht]
\caption{Average GMRES iterations and value of the discrete cost functional $J^*_h$ with the augmented Lagrangian preconditioner with  $\gamma=10/\sqrt{\beta}$, for $\nu=\frac{1}{250}$ and a range of $l$, $\beta$.}\label{Inexact_Newton_Opt_then_Discr_Augmented_Lagrangian_2}
\begin{scriptsize}
\begin{center}
\renewcommand{\arraystretch}{1.2}
{\begin{tabular}{|c||c|c||c|c||c|c||c|c||c|c||c|c||c|c|}
\hline
\multicolumn{1}{|c||}{} & \multicolumn{14}{c|}{$\beta$} \\
\cline{2-15}
\multicolumn{1}{|c||}{} & \multicolumn{2}{c||}{$10^{-1}$} & \multicolumn{2}{c||}{$10^{-2}$} & \multicolumn{2}{c||}{$10^{-3}$} & \multicolumn{2}{c||}{$10^{-4}$} & \multicolumn{2}{c||}{$10^{-5}$} & \multicolumn{2}{c||}{$10^{-6}$} & \multicolumn{2}{c|}{$10^{-8}$} \\
\cline{2-15}
$l$ & $\mathrm{it}$ & $J^*_h$ & $\mathrm{it}$ & $J^*_h$ & $\mathrm{it}$ & $J^*_h$ & $\mathrm{it}$ & $J^*_h$ & $\mathrm{it}$ & $J^*_h$ & $\mathrm{it}$ & $J^*_h$ & $\mathrm{it}$ & $J^*_h$\\
\hline
\hline
\!\!\!\! $3$ \!\!\!\! & 6 & 3.5e-02 & 5 & 3.1e-02 & 4 & 3.1e-02 & 4 & 3.0e-02 & 3 & 3.0e-02 & 3 & 3.0e-02 & 2 & 3.0e-02 \\
\hline
\!\!\!\! $4$ \!\!\!\! & 7 & 2.4e-02 & 6 & 1.8e-02 & 4 & 1.5e-02 & 4 & 1.5e-02 & 3 & 1.5e-02 & 3 & 1.5e-02 & 3 & 1.5e-02 \\
\hline
\!\!\!\! $5$ \!\!\!\! & 7 & 2.4e-02 & 7 & 1.4e-02 & 5 & 9.3e-03 & 5 & 7.8e-03 & 3 & 7.5e-03 & 3 & 7.5e-03 & 3 & 7.5e-03 \\
\hline
\!\!\!\! $6$ \!\!\!\! & 6 & 2.4e-02 & 6 & 1.4e-02 & 6 & 7.8e-03 & 5 & 4.9e-03 & 5 & 3.9e-03 & 4 & 3.7e-03 & 4 & 3.7e-03 \\
\hline
\!\!\!\! $7$ \!\!\!\! & 7 & 2.4e-02 & 6 & 1.4e-02 & 6 & 7.8e-03 & 6 & 4.4e-03 & 5 & 2.7e-03 & 5 & 2.0e-03 & 4 & 1.9e-03 \\
\hline
\end{tabular}}
\end{center}
\end{scriptsize}
\end{table}

\begin{table}[!ht]
\caption{Average GMRES iterations and value of the discrete cost functional $J^*_h$ with the augmented Lagrangian preconditioner with  $\gamma=10/\sqrt{\beta}$, for $\nu=\frac{1}{500}$ and a range of $l$, $\beta$.}\label{Inexact_Newton_Opt_then_Discr_Augmented_Lagrangian_3}
\begin{scriptsize}
\begin{center}
\renewcommand{\arraystretch}{1.2}
{\begin{tabular}{|c||c|c||c|c||c|c||c|c||c|c||c|c||c|c|}
\hline
\multicolumn{1}{|c||}{} & \multicolumn{14}{c|}{$\beta$} \\
\cline{2-15}
\multicolumn{1}{|c||}{} & \multicolumn{2}{c||}{$10^{-1}$} & \multicolumn{2}{c||}{$10^{-2}$} & \multicolumn{2}{c||}{$10^{-3}$} & \multicolumn{2}{c||}{$10^{-4}$} & \multicolumn{2}{c||}{$10^{-5}$} & \multicolumn{2}{c||}{$10^{-6}$} & \multicolumn{2}{c|}{$10^{-8}$} \\
\cline{2-15}
$l$ & $\mathrm{it}$ & $J^*_h$ & $\mathrm{it}$ & $J^*_h$ & $\mathrm{it}$ & $J^*_h$ & $\mathrm{it}$ & $J^*_h$ & $\mathrm{it}$ & $J^*_h$ & $\mathrm{it}$ & $J^*_h$ & $\mathrm{it}$ & $J^*_h$\\
\hline
\hline
\!\!\!\! $3$ \!\!\!\! & 7 & 3.2e-02 & 5 & 3.1e-02 & 4 & 3.0e-02 & 4 & 3.0e-02 & 3 & 3.0e-02 & 3 & 3.0e-02 & 2 & 3.0e-02 \\
\hline
\!\!\!\! $4$ \!\!\!\! & 9 & 1.9e-02 & 6 & 1.6e-02 & 4 & 1.5e-02 & 4 & 1.5e-02 & 3 & 1.5e-02 & 3 & 1.5e-02 & 3 & 1.5e-02 \\
\hline
\!\!\!\! $5$ \!\!\!\! & 8 & 1.7e-02 & 6 & 1.0e-02 & 5 & 8.1e-03 & 5 & 7.6e-03 & 3 & 7.5e-03 & 3 & 7.5e-03 & 3 & 7.5e-03 \\
\hline
\!\!\!\! $6$ \!\!\!\! & 8 & 1.7e-02 & 6 & 9.7e-03 & 5 & 5.7e-03 & 5 & 4.2e-03 & 5 & 3.8e-03 & 4 & 3.7e-03 & 4 & 3.7e-03 \\
\hline
\!\!\!\! $7$ \!\!\!\! & 9 & 1.7e-02 & 6 & 9.8e-03 & 6 & 5.5e-03 & 7 & 3.2e-03 & 5 & 2.2e-03 & 5 & 1.9e-03 & 4 & 1.9e-03 \\
\hline
\end{tabular}}
\end{center}
\end{scriptsize}
\end{table}

\begin{table}[!ht]
\caption{Degrees of freedom (DoF) and number of inexact Newton iterations required for stationary Navier--Stokes control problem. In each cell are the inexact Newton iterations for the given $l$, $\nu$, and $\beta=10^{-j}$, $j=1,\ldots ,5$.}\label{Number_Inexact_newton_Opt_then_Discr}
\begin{scriptsize}
\begin{center}
\renewcommand{\arraystretch}{1.2}
\begin{tabular}{|c|c||ccccc|ccccc|ccccc|}
\hline
$l$ & DoF & \multicolumn{5}{c|}{$\nu=\frac{1}{100}$} & \multicolumn{5}{c|}{$\nu=\frac{1}{250}$} & \multicolumn{5}{c|}{$\nu=\frac{1}{500}$} \\
\hline
\hline
$3$ & 1062 & 5 & 4 & 3 & 3 & 3 & 7 & 5 & 4 & 3 & 3 & 8 & 5 & 4 & 3 & 3 \\
\hline
$4$ & 4422 & 4 & 4 & 3 & 3 & 3 & 6 & 4 & 3 & 3 & 3 & 8 & 5 & 5 & 4 & 3 \\
\hline
$5$ & 18,054 & 4 & 3 & 3 & 3 & 3 & 4 & 4 & 3 & 3 & 3 & 6 & 4 & 3 & 3 & 3 \\
\hline
$6$ & 72,966 & 4 & 3 & 3 & 4 & 3 & 4 & 3 & 3 & 3 & 3 & 4 & 3 & 3 & 3 & 3 \\
\hline
$7$ & 293,382 & 3 & 3 & 3 & 3 & 3 & 3 & 3 & 3 & 3 & 3 & 4 & 3 & 3 & 3 & 3 \\
\hline
\end{tabular}
\end{center}
\end{scriptsize}
\end{table}

\subsubsection{Discretize-then-Optimize}\label{sec_4_1_2}
In this section, we report the results of our augmented Lagrangian
	preconditioner when employing an inexact solver on the augmented
	block, when considering a discretize-then-optimize approach.
	We employ the Firedrake system
	\cite{Rathgeber_Ham_Mitchell_Lange_Luporini_Mcrae_Bercea_Markall_Kelley}
	for the construction of the finite element discretizations of the
	inexact Newton iteration and of the linear and non-linear solvers employed.
	We employ the exactly incompressible Scott--Vogelius $[P_4]^2$--$DP_3$
	finite element pair~\cite{Scott_Vogelius}. We mention that for this problem
	no stabilization for the linearized convection term is employed in the
	discretization. Since we employ
	discontinuous finite elements on the pressure space, the mass matrix is
	block-diagonal and is solved cellwise. The stiffness matrix is solved with a
	two-level $p$-multigrid solver applied to an auxiliary-space
	$P_3$ discretization, i.e.~$DP_3 \rightarrow P_3 \rightarrow P_1$, with
	Chebyshev-accelerated symmetric Gauss--Seidel pointwise relaxation on the
	degree-3 levels and a geometric multigrid cycle applied on the $P_1$ problem.
	Finally, the momentum blocks are inexactly solved with one cycle of the
	multigrid derived in \cite{Benzi_Olshanskii_2006} and implemented in
	Firedrake in \cite{Farrell_Mitchell_Wechsung}. Ten iterations of FGMRES
	preconditioned by an additive vertex-star solve are employed as relaxation
	on each fine level, while the coarse problem is solved with
	SuperLU\_DIST~\cite{Li_Demmel}. The vertex-star solve is performed using
	PETSc's PCASM~\cite{PETSc,Smith_Bjorstad_Gropp}. For the coarsest
	level ($l=1$), multigrid has been replaced by an exact solver based on the
	LU factorization of each block. We employ a Intel Core Ultra 7 165U
	processor and 64 GB RAM on an Ubuntu 22.04.5 LTS operating system.

We allow for 10 inexact Newton iterations, specifying as a stopping
	criterion a reduction of $10^{-5}$ on either the relative or the
	absolute non-linear residual. We run our
	preconditioned iterative solver up to a tolerance of $10^{-6}$ on
	either the relative or the absolute residual is achieved. We report
	in Tables \ref{Inexact_Newton_Discr_then_Opt_1}--\ref{Inexact_Newton_Discr_then_Opt_2}
	the average number of FGMRES iterations required for reaching
	convergence, together with the value of the cost
	functional $J^*_h$. Further, we report in Tables
	\ref{Inexact_Newton_Discr_then_Opt_CPU_1}--\ref{Inexact_Newton_Discr_then_Opt_CPU_3} the total CPU times required for the
	solution of the problem, with the dimension of the systems solved showed
	and of the average block-smoother per multigrid level
	in Table \ref{Inexact_Newton_Discr_then_Opt_CPU_1}. Note that the CPU
	times reported include the evaluation of the non-linear residual, the
	discretization of the operators, and the construction of the geometric
	multigrid on the augmented velocity blocks, in addition to the times
	required for the linear solvers. Further, note that $l=1$ represents the coarsest mesh, thus no smoother has been employed.
	
From Tables \ref{Inexact_Newton_Discr_then_Opt_1}--\ref{Inexact_Newton_Discr_then_Opt_2},
	we can observe again the overall
	robustness of the augmented Lagrangian approach, even when employing
	inexact solvers throughout.
	Only for larger values of $\beta$ on the finest grid do we
	see some growth in the number of FGMRES iterations; nevertheless,
	even in this setting the solver
	requires at most 15 iterations on average to reach convergence. Further, we observe that the discrete cost
	functional $J^*_h$ tends to a minimum as the mesh is refined. Finally,
	from Tables \ref{Inexact_Newton_Discr_then_Opt_CPU_1}--\ref{Inexact_Newton_Discr_then_Opt_CPU_3} 
	we observe that for the finest meshes the CPU time scales linearly with
	respect to the problem size, while for coarsest meshes we observe the smal
	problem scaling effect.

	\begin{table}[!ht]
		\caption{Average FGMRES iterations and value of the discrete cost
		functional $J^*_h$ with the augmented Lagrangian
		preconditioner with $\gamma=10/\sqrt{\beta}$, for
		$\nu=\frac{1}{100}$, and a
		range of $l$, $\beta$.}\label{Inexact_Newton_Discr_then_Opt_1}
		
		\begin{scriptsize}
			\begin{center}
				\renewcommand{\arraystretch}{1.2}
				{\begin{tabular}{|c||c|c||c|c||c|c||c|c||c|c||c|c||c|c|}
					\hline
					 & \multicolumn{14}{c|}{$\beta$} \\
					\cline{2-15}
					 & \multicolumn{2}{c||}{$10^{-1}$} & \multicolumn{2}{c||}{$10^{-2}$} & \multicolumn{2}{c||}{$10^{-3}$} & \multicolumn{2}{c||}{$10^{-4}$} & \multicolumn{2}{c||}{$10^{-5}$} & \multicolumn{2}{c||}{$10^{-6}$} & \multicolumn{2}{c|}{$10^{-8}$} \\
					\cline{2-15}
					$l$ & $\mathrm{it}$ & $J^*_h$ & $\mathrm{it}$ & $J^*_h$ & $\mathrm{it}$ & $J^*_h$ & $\mathrm{it}$ & $J^*_h$ & $\mathrm{it}$ & $J^*_h$ & $\mathrm{it}$ & $J^*_h$ & $\mathrm{it}$ & $J^*_h$\\
					\hline\hline
					$1$ & 3 & 4.5e-02 & 3 & 4.3e-02 & 3 & 4.3e-02 & 3 & 4.3e-02 & 3 & 4.3e-02 & 3 & 4.3e-02 & 2 & 4.2e-02 \\
					\hline
					$2$ & 3 & 3.2e-02 & 3 & 2.3e-02 & 3 & 2.1e-02 & 2 & 2.1e-02 & 3 & 2.1e-02 & 3 & 2.1e-02 & 2 & 2.1e-02 \\
					\hline
					$3$ & 3 & 3.4e-02 & 3 & 2.0e-02 & 3 & 1.2e-02 & 4 & 1.1e-02 & 3 & 0.1e-02 & 3 & 1.0e-02 & 2 & 1.0e-02 \\
					\hline
					$4$ & 5 & 3.5e-02 & 4 & 2.1e-02 & 3 & 1.1e-02 & 3 & 6.8e-03 & 4 & 5.4e-03 & 3 & 5.2e-03 & 2 & 5.1e-03 \\
					\hline
					$5$ & 6 & 3.6e-02 & 7 & 2.1e-02 & 3 & 1.2e-02 & 3 & 6.8e-03 & 3 & 3.8e-03 & 3 & 2.8e-03 & 3 & 2.6e-03 \\
					\hline
					$6$ & 9 & 3.6e-02 & 5 & 2.1e-02 & 3 & 1.2e-02 & 3 & 6.9e-03 & 3 & 3.9e-03 & 3 & 2.1e-03 & 4 & 1.3e-03 \\
					\hline
				\end{tabular}}
			\end{center}
		\end{scriptsize}
	\end{table}

	\begin{table}[!ht]
		\caption{Average FGMRES iterations and value of the discrete cost
		functional $J^*_h$ with the augmented Lagrangian
		preconditioner with $\gamma=10/\sqrt{\beta}$, for
		$\nu=\frac{1}{250}$, and a
		range of $l$, $\beta$.}\label{Inexact_Newton_Discr_then_Opt_2}
		
		\begin{scriptsize}
			\begin{center}
				\renewcommand{\arraystretch}{1.2}
				{\begin{tabular}{|c||c|c||c|c||c|c||c|c||c|c||c|c||c|c|}
					\hline
					\multicolumn{1}{|c||}{} & \multicolumn{14}{c|}{$\beta$} \\
					\cline{2-15}
					 & \multicolumn{2}{c||}{$10^{-1}$} & \multicolumn{2}{c||}{$10^{-2}$} & \multicolumn{2}{c||}{$10^{-3}$} & \multicolumn{2}{c||}{$10^{-4}$} & \multicolumn{2}{c||}{$10^{-5}$} & \multicolumn{2}{c||}{$10^{-6}$} & \multicolumn{2}{c|}{$10^{-8}$} \\
					\cline{2-15}
					$l$ & $\mathrm{it}$ & $J^*_h$ & $\mathrm{it}$ & $J^*_h$ & $\mathrm{it}$ & $J^*_h$ & $\mathrm{it}$ & $J^*_h$ & $\mathrm{it}$ & $J^*_h$ & $\mathrm{it}$ & $J^*_h$ & $\mathrm{it}$ & $J^*_h$\\
					\hline
					\hline
					$1$ & 3 & 4.4e-02 & 3 & 4.3e-02 & 3 & 4.3e-02 & 3 & 4.3e-02 & 3 & 4.3e-02 & 3  & 4.3e-02 & 2 & 4.2e-02 \\
					\hline
					$2$ & 3 & 2.4e-02 & 3 & 2.2e-02 & 3 & 2.1e-02 & 3 & 2.1e-02 & 3 & 2.1e-02 & 3 & 2.1e-02 & 2 & 2.1e-02 \\
					\hline
					$3$ & 4 & 2.2e-02 & 3 & 1.3e-02 & 4 & 1.1e-02 & 3 & 1.0e-02 & 3 & 1.0e-02 & 2 & 1.0e-02 & 2 & 1.0e-02 \\
					\hline
					$4$ & 4 & 2.3e-02 & 4 & 1.3e-02 & 3 & 7.5e-03 & 3 & 5.5e-03 & 4 & 5.2e-03 & 4 & 5.1e-03 & 1 & 5.0e-03 \\
					\hline
					$5$ & 6 & 2.3e-02 & 5 & 1.3e-02 & 3 & 7.6e-03 & 3 & 4.2e-03 & 3 & 2.9e-03 & 4 & 2.6e-03 & 4 & 2.6e-03 \\
					\hline
					$6$ & 13 & 2.3e-02 & 8 & 1.4e-02 & 3 & 7.7e-03 & 3 & 4.4e-03 & 3 & 2.3e-03 & 4 & 1.5e-03 & 4 & 1.3e-03 \\
					\hline
				\end{tabular}}
			\end{center}
		\end{scriptsize}
	\end{table}

	\begin{table}[!ht]
		\caption{Average FGMRES iterations and value of the discrete cost
		functional $J^*_h$ with the augmented Lagrangian preconditioner with
		$\gamma=10/\sqrt{\beta}$, for $\nu=\frac{1}{500}$, and a
		range of $l$, $\beta$.}\label{Inexact_Newton_Discr_then_Opt_3}
		
		\begin{scriptsize}
			\begin{center}
				\renewcommand{\arraystretch}{1.2}
				\begin{tabular}{|c||c|c||c|c||c|c||c|c||c|c||c|c||c|c|}
					\hline
					\multicolumn{1}{|c||}{} & \multicolumn{14}{c|}{$\beta$} \\
					\cline{2-15}
					 & \multicolumn{2}{c||}{$10^{-1}$} & \multicolumn{2}{c||}{$10^{-2}$} & \multicolumn{2}{c||}{$10^{-3}$} & \multicolumn{2}{c||}{$10^{-4}$} & \multicolumn{2}{c||}{$10^{-5}$} & \multicolumn{2}{c||}{$10^{-6}$} & \multicolumn{2}{c|}{$10^{-8}$} \\
					\cline{2-15}
					$l$ & $\mathrm{it}$ & $J^*_h$ & $\mathrm{it}$ & $J^*_h$ & $\mathrm{it}$ & $J^*_h$ & $\mathrm{it}$ & $J^*_h$ & $\mathrm{it}$ & $J^*_h$ & $\mathrm{it}$ & $J^*_h$ & $\mathrm{it}$ & $J^*_h$\\
					\hline\hline
					$1$ & 3 & 4.3e-02 & 3 & 4.3e-02 & 3 & 4.3e-02 & 3 & 4.3e-02 & 3 & 4.3e-02 & 3 & 4.3e-02 & 2 & 4.2e-02 \\
					\hline
					$2$ & 3 & 2.2e-02 & 3 & 2.1e-02 & 2 & 2.1e-02 & 3 & 2.1e-02 & 3 & 2.1e-02 & 3 & 2.1e-02 & 2 & 2.1e-02 \\
					\hline
					$3$ & 4 & 1.6e-02 & 4 & 1.1e-02 & 4 & 1.0e-02 & 3 & 1.0e-02 & 3 & 1.0e-02 & 2 & 1.0e-02 & 2 & 1.0e-02 \\
					\hline
					$4$ & 4 & 1.6e-02 & 3 & 9.2e-03 & 3 & 6.0e-03 & 4 & 5.2e-03 & 4 & 5.2e-03 & 4 & 5.1e-03 & 4 & 5.1e-03 \\
					\hline
					$5$ & 8 & 1.7e-02 & 5 & 9.6e-03 & 3 & 5.3e-03 & 3 & 3.2e-03 & 4 & 2.7e-03 & 4 & 2.6e-03 & 2 & 2.5e-03 \\
					\hline
					$6$ & 15 & 1.7e-02 & 7 & 9.7e-03 & 3 & 5.5e-03 & 3 & 3.1e-03 & 3 & 1.7e-03 & 4 & 1.4e-03 & 4 & 1.3e-03 \\
					\hline
				\end{tabular}
			\end{center}
		\end{scriptsize}
	\end{table}

	\begin{table}[!ht]
	\caption{Degrees of freedom (DoF), average size of the block-smoothers (BlkSm), and total CPU times, for the augmented Lagrangian preconditioner with $\gamma=10/\sqrt{\beta}$, for $\nu=\frac{1}{100}$ and a range of $l$, $\beta$.}\label{Inexact_Newton_Discr_then_Opt_CPU_1}
		\begin{scriptsize}
			\begin{center}
			\renewcommand{\arraystretch}{1.3}
				{\begin{tabular}{|c||c||c||c|c|c|c|c|c|c|}
					\hline
					 & & & \multicolumn{7}{c|}{$\beta$} \\
					\cline{4-10}
					\!\!\!\! $l$ \!\!\!\! & DoF & BlkSm & $10^{-1}$ & $10^{-2}$ & $10^{-3}$ & $10^{-4}$ & $10^{-5}$ & $10^{-6}$ & $10^{-8}$ \\
					\hline
					\hline
					$1$ & 484 & \textbackslash & 16.7 & 15.6 & 16.1 & 13.4 & 12.8 & 12.5 & 12.3 \\
					\hline
					$2$ & 1796 & 61 & 25.2 & 26.1 & 26.6 & 21.1 & 19.6 & 22.1 & 19.4 \\
					\hline
					$3$ & 6916 & 67 & 34.8 & 27.7 & 20.5 & 21.6 & 21.8 & 22.5 & 19.9 \\
					\hline
					$4$ & 27,140 & 70 & 35.7 & 33.9 & 24.4 & 24.6 & 25.6 & 24.3 & 23.4 \\
					\hline
					$5$ & 107,524 & 72 & 61.5 & 67.9 & 33.1 & 32.8 & 33.5 & 33.9 & 32.5 \\
					\hline
					$6$ & 428,036 & 73 & 218 & 136 & 61.6 & 61.9 & 63.3 & 62.8 & 73.1 \\
					\hline
				\end{tabular}}
			\end{center}
		\end{scriptsize}
	\end{table}

\begin{table}[!ht]
\caption{Total CPU times, for the augmented Lagrangian preconditioner with $\gamma=10/\sqrt{\beta}$, for $\nu=\frac{1}{250}$ and a range of $l$, $\beta$.}\label{Inexact_Newton_Discr_then_Opt_CPU_2}
\begin{scriptsize}
\begin{center}
\renewcommand{\arraystretch}{1.3}
{\begin{tabular}{|c||c|c|c|c|c|c|c|}
\hline
\multicolumn{1}{|c||}{} & \multicolumn{7}{c|}{$\beta$} \\
\cline{2-8}
$l$ & $10^{-1}$ & $10^{-2}$ & $10^{-3}$ & $10^{-4}$ & $10^{-5}$ & $10^{-6}$ & $10^{-8}$ \\
\hline
\hline
$1$ & 16.2 & 16.5 & 16.4 & 12.9 & 12.9 & 13.2 & 12.9  \\
\hline
$2$ & 30.6 & 26.8 & 26.6 & 19.6 & 20.5 & 18.9 & 19.7 \\
\hline
$3$ & 28.5 & 21.5 & 23.4 & 21.2 & 21.7 & 20.8 & 21.3 \\
\hline
$4$ & 47.1 & 35.1 & 24.1 & 25.1 & 26.3 & 25.5 & 22.2 \\
\hline
$5$ & 65.8 & 56.8 & 32.4 & 33.3 & 34.3 & 36.1 & 36.8 \\
\hline
$6$ & 294 & 195 & 61.6 & 62.7 & 61.9 & 73.5 & 73.6 \\
\hline
\end{tabular}}
\end{center}
\end{scriptsize}
\end{table}

\begin{table}[!ht]
\caption{Total CPU times, for the augmented Lagrangian preconditioner with $\gamma=10/\sqrt{\beta}$, for $\nu=\frac{1}{500}$ and a range of $l$, $\beta$.}\label{Inexact_Newton_Discr_then_Opt_CPU_3}
\begin{scriptsize}
\begin{center}
\renewcommand{\arraystretch}{1.3}
{\begin{tabular}{|c||c|c|c|c|c|c|c|}
\hline
\multicolumn{1}{|c||}{} & \multicolumn{7}{c|}{$\beta$} \\
\cline{2-8}
$l$ & $10^{-1}$ & $10^{-2}$ & $10^{-3}$ & $10^{-4}$ & $10^{-5}$ & $10^{-6}$ & $10^{-8}$ \\
\hline
\hline
$1$ & 15.9 & 16.3 & 16.7 & 12.7 & 12.8 & 12.9 & 12.7 \\
\hline
$2$ & 26.4 & 28.1 & 26.9 & 19.4 & 19.8 & 19.7 & 19.9 \\
\hline
$3$ & 30.4 & 30.7 & 22.7 & 21.9 & 21.3 & 21.5 & 22.4 \\
\hline
$4$ & 45.5 & 25.4 & 26.1 & 26.3 & 26.4 & 26.2 & 25.1 \\
\hline
$5$ & 76.2 & 57.9 & 33.5 & 33.7 & 38.3 & 36.5 & 31.3 \\
\hline
$6$ & 338 & 185 & 62.3 & 64.7 & 63.4 & 74.2 & 74.3 \\
\hline
\end{tabular}}
\end{center}
\end{scriptsize}
\end{table}

\subsection{Backward-Facing Step}\label{sec_4_2}
As last test, we solve the Navier--Stokes control problem for the backward-facing
	step. We apply a Poiseuille inflow on the left end and ``no-slip''
	(zero boundary) conditions on the top and bottom of the step. Natural boundary
	conditions are imposed on the right end of the step. We set $\vec{f}=\vec{0}$,
	and seek as desired state $\vec{v}_d$ the solution of the stationary
	Stokes problem obtained when imposing this forcing and boundary conditions.
	
As in Section \ref{sec_4_1_2}, we apply the discretize-then-optimize strategy,
	employing the Firedrake system for the derivation of the optimality conditions
	of the linearized problem and for the construction of the linear and non-linear
	solvers. Further, as above we employ the Scott--Vogelius $[P_4]^2$--$DP_3$
	finite element pair. Finally, we allow for 10 inexact Newton iterations,
	specifying as a stopping criterion a reduction of $10^{-5}$ on either the
	relative or the absolute non-linear residual, running the linear solver up to
	a tolerance $10^{-6}$ on either the relative or the absolute linear residual.
	For this test, we set $\gamma=10^3$ and $\gamma=10^4$. We opt for this choice
	of $\gamma$ as the step problem is more difficult to solve, due to
	the singularity of the pressure in the corner of the step. We report in
	Table \ref{step_gamma1000} the results for
	$\gamma=10^3$ and in Table \ref{step_gamma10000}
	the results for $\gamma=10^4$, respectively. Further, in Tables
	\ref{CPU_step_gamma1000}--\ref{CPU_step_gamma10000} we report the total CPU
	times required for solving the problem, with the dimension of each linear
	system solved and of the average block-smoother per multigrid level reported
	in Table \ref{CPU_step_gamma1000}.

\begin{table}[!ht]
\caption{Inexact Newton iteration: average FGMRES iterations and value of the discrete cost functional $J^*_h$ with the augmented Lagrangian preconditioner with $\gamma=10^3$, for $\nu=\frac{1}{100}$, $\frac{1}{250}$, and $\frac{1}{500}$, and a range of $l$, $\beta$.}\label{step_gamma1000}

\begin{scriptsize}
\begin{center}
\renewcommand{\arraystretch}{1.3}
{\begin{tabular}{|c||c||c|c||c|c||c|c||c|c||c|c||c|c||c|c|}
\hline
\multicolumn{1}{|c||}{} & \multicolumn{1}{c||}{} & \multicolumn{14}{c|}{$\beta$} \\
\cline{3-16}
 & & \multicolumn{2}{c||}{$10^{-1}$} & \multicolumn{2}{c||}{$10^{-2}$} & \multicolumn{2}{c||}{$10^{-3}$} & \multicolumn{2}{c||}{$10^{-4}$} & \multicolumn{2}{c||}{$10^{-5}$} & \multicolumn{2}{c||}{$10^{-6}$} & \multicolumn{2}{c|}{$10^{-8}$} \\
\cline{3-16}
$\nu$ & $l$ & $\mathrm{it}$ & $J^*_h$ & $\mathrm{it}$ & $J^*_h$ & $\mathrm{it}$ & $J^*_h$ & $\mathrm{it}$ & $J^*_h$ & $\mathrm{it}$ & $J^*_h$ & $\mathrm{it}$ & $J^*_h$ & $\mathrm{it}$ & $J^*_h$ \\
\hline\hline
\multirow{3}{*}{$\frac{1}{100}$} & 
$1$ & 3 & 1.3e00 & 3 & 1.3e00 & 3 & 2.6e00 & 3 & 1.2e00 & 5 & 6.9e-01 & 5 & 5.3e-01 & 13 & 3.0e-01 \\
\cline{3-16}
& $2$ & 4 & 1.4e00 & 3 & 1.3e00 & 5 & 1.6e00 &4 & 1.4e00 & 8 & 8.4e-01 & 20 & 1.6e00 & 12 & 3.2e-01 \\
\cline{3-16}
& $3$ & 3 & 1.6e00 & 3 & 1.5e00 & 3 & 2.5e00 & 3 & 1.2e00 & 6 & 7.1e-01 & 11 & 6.0e-01 & 12 & 2.7e-01 \\
\hline
\multirow{3}{*}{$\frac{1}{250}$} & 
$1$ & 3 & 1.3e00 & 3 & 1.3e00 & 3 & 2.6e00 & 3 & 1.2e00 & 5 & 6.9e-01 & 5 & 5.3e-01 & 13 & 3.0e-01 \\
\cline{3-16}
& $2$ & 4 & 1.4e00 & 3 & 1.4e00 & 5 & 1.6e00 & 10 & 1.4e00 & 16 & 1.9e00 & 22 & 2.2e00 & 12 & 3.2e-01 \\
\cline{3-16}
& $3$ &3 & 1.7e00 & 3 & 1.5e00 & 4 & 2.3e00 & 5 & 1.8e00 & 9 & 9.6e-01 & 20 & 8.9e-01 & 12 & 3.2e-01 \\
\hline
\multirow{3}{*}{$\frac{1}{500}$} & 
$1$ & 3 & 1.3e00 & 3 & 1.3e00 & 3 & 2.6e00 & 3 & 1.2e00 & 5 & 6.9e-01 & 5 & 5.3e-01 & 13 & 3.0e-01 \\
\cline{3-16}
& $2$ & 3 & 1.7e00 & 4 & 1.3e00 & 6 & 1.4e00 & 9 & 1.5e00 & 17 & 1.7e00 & 17 & 6.1e-01 & 12 & 3.1e-01 \\
\cline{3-16}
& $3$ & 3 & 1.7e00 & 3 & 1.5e00 & 5 & 1.8e00 & 8 & 1.8e00 & 14 & 1.1e00 & 22 & 2.4e00 & 11 & 3.2e-01 \\
\hline
\end{tabular}}
\end{center}
\end{scriptsize}
\end{table}

\begin{table}[!ht]
\caption{Inexact Newton iteration: average FGMRES iterations and value of the discrete cost functional $J^*_h$ with the augmented Lagrangian preconditioner with $\gamma=10^4$, for $\nu=\frac{1}{100}$, $\frac{1}{250}$, and $\frac{1}{500}$, and a range of $l$, $\beta$.}\label{step_gamma10000}

\begin{scriptsize}
\begin{center}
\renewcommand{\arraystretch}{1.3}
{\begin{tabular}{|c||c||c|c||c|c||c|c||c|c||c|c||c|c||c|c|}
\hline
\multicolumn{1}{|c||}{} & \multicolumn{1}{c||}{} & \multicolumn{14}{c|}{$\beta$} \\
\cline{3-16}
 & & \multicolumn{2}{c||}{$10^{-1}$} & \multicolumn{2}{c||}{$10^{-2}$} & \multicolumn{2}{c||}{$10^{-3}$} & \multicolumn{2}{c||}{$10^{-4}$} & \multicolumn{2}{c||}{$10^{-5}$} & \multicolumn{2}{c||}{$10^{-6}$} & \multicolumn{2}{c|}{$10^{-8}$} \\
\cline{3-16}
$\nu$ & $l$ & $\mathrm{it}$ & $J^*_h$ & $\mathrm{it}$ & $J^*_h$ & $\mathrm{it}$ & $J^*_h$ & $\mathrm{it}$ & $J^*_h$ & $\mathrm{it}$ & $J^*_h$ & $\mathrm{it}$ & $J^*_h$ & $\mathrm{it}$ & $J^*_h$ \\
\hline\hline
\multirow{3}{*}{$\frac{1}{100}$} & 
$1$ & 2 & 3.5e00 & 3 & 1.3e00 & 2 & 3.3e00 & 2 & 3.2e00 & 3 & 2.5e00 & 2 & 1.3e00 & 4 & 5.3e-01 \\
\cline{3-16}
& $2$ & 2 & 2.4e00 & 2 & 2.4e00 & 3 & 1.7e00 & 4 & 1.5e00 & 7 & 1.2e00 & 10 & 1.4e00 & 22 & 3.4e00 \\
\cline{3-16}
& $3$ & 3 & 1.8e00 & 3 & 1.7e00 & 3 & 1.7e00 & 3 & 1.3e00 & 5 & 1.7e00 & 10 & 1.3e00 & 22 & 3.1e00 \\
\hline
\multirow{3}{*}{$\frac{1}{250}$} & 
$1$ & 2 & 3.4e00 & 3 & 1.3e00 & 2 & 3.3e00 & 2 & 3.2e00 & 3 & 2.5e00 & 2 & 1.3e00 & 4 & 5.2e-01 \\
\cline{3-16}
& $2$ & 2 & 2.3e00 & 2 & 2.3e00 & 4 & 1.3e00 & 4 & 1.3e00 & 8 & 1.3e00 & 10 & 1.4e00 & 21 & 2.7e00 \\
\cline{3-16}
& $3$ & 2 & 2.4e00 & 3 & 1.7e00 & 3 & 1.3e00 & 4 & 1.5e00 & 7 & 1.3e00 & 14 & 1.3e00 & 16 & 6.1e-01 \\
\hline
\multirow{3}{*}{$\frac{1}{500}$} & 
$1$ & 2 & 3.4e00 & 3 & 1.3e00 & 2 & 3.3e00 & 2 & 3.2e00 & 3 & 2.5e00 & 2 & 1.3e00 &4 & 5.2e-01 \\
\cline{3-16}
& $2$ & 3 & 1.7e00 & 3 & 1.7e00 & 3 & 1.7e00 & 6 & 1.3e00 & 7 & 1.3e00 & 7 & 2.1e00 & 22 & 1.9e00 \\
\cline{3-16}
& $3$ & 3 & 1.7e00 & 3 & 1.8e00 & 3 & 1.7e00 & 5 & 1.3e00 & 7 & 1.2e00 & 10 & 1.4e00 & 19 & 1.9e00 \\
\hline
\end{tabular}}
\end{center}
\end{scriptsize}
\end{table}

From Tables \ref{step_gamma1000}--\ref{step_gamma10000}, we observe satisfactory
	robustness of the solver for the larger values of $\beta$, with a small increase
	in the number of iterations for $\beta=10^{-6}$ and $\beta=10^{-8}$.
	Nonetheless, the solver is able to obtain convergence in at most 22 iterations
	on average. Further, from Tables
	\ref{CPU_step_gamma1000}--\ref{CPU_step_gamma10000} we observe that the
	CPU times scale nearly linearly with respect to the problem size.

\begin{table}[!ht]
\caption{Inexact Newton iteration: degrees of freedom (DoF), average size of the block-smoothers (BlkSm), and total CPU times of the augmented Lagrangian preconditioner with $\gamma=10^3$, for $\nu=\frac{1}{100}$, $\frac{1}{250}$, and $\frac{1}{500}$, and a range of $l$, $\beta$.}\label{CPU_step_gamma1000}

\begin{scriptsize}
\begin{center}
\renewcommand{\arraystretch}{1.2}
{\begin{tabular}{|c||c||c||c||c|c|c|c|c|c|c|}
\hline
\multicolumn{1}{|c||}{} & \multicolumn{1}{|c||}{} & & \multicolumn{1}{c||}{} & \multicolumn{7}{c|}{$\beta$} \\
\cline{5-11}
$\nu$ & $l$ & DoF & BlkSm & $10^{-1}$ & $10^{-2}$ & $10^{-3}$ & $10^{-4}$ & $10^{-5}$ & $10^{-6}$ & $10^{-8}$ \\
\hline\hline
\multirow{3}{*}{$\frac{1}{100}$} & 
$1$ & 168,184 & \textbackslash & 11.5 & 3.8 & 3.7 & 4.0 & 4.5 & 4.7 & 7.1 \\
\cline{4-11}
& $2$ & 669,476 & 72 & 81.0 & 58.9 & 81.1 & 70.2 & 116 & 258 & 164 \\
\cline{4-11}
& $3$ & 2,671,396 & 73 & 364 & 358 & 360 & 358 & 555 & 884 & 951 \\
\hline
\multirow{3}{*}{$\frac{1}{250}$} & 
$1$ & 168,184 & \textbackslash & 3.6 & 3.7 & 3.7 & 4.0 & 4.6 & 4.8 & 7.2 \\
\cline{4-11}
& $2$ & 669,476 & 72 & 69.5 & 58.5 & 81.7 & 140 & 212 & 283 & 165 \\
\cline{4-11}
& $3$ & 2,671,396 & 73 & 360 & 359 & 426 & 492 & 757 & 1481 & 948 \\
\hline
\multirow{3}{*}{$\frac{1}{500}$} & 
$1$ & 168,184 & \textbackslash & 3.5 & 3.8 & 3.8 & 3.9 & 4.8 & 4.6 & 7.2 \\
\cline{4-11}
& $2$ & 669,476 & 72 & 58.8 & 70.7 & 93.6 & 129 & 224 & 224 & 166 \\
\cline{4-11}
& $3$ & 2,671,396 & 73 & 361 & 359 & 491 & 684 & 1080 & 1604 & 884 \\
\hline
\end{tabular}}
\end{center}
\end{scriptsize}
\end{table}

\begin{table}[!ht]
\caption{Inexact Newton iteration: total CPU times of the augmented Lagrangian preconditioner with $\gamma=10^4$, for $\nu=\frac{1}{100}$, $\frac{1}{250}$, and $\frac{1}{500}$, and a range of $l$, $\beta$.}\label{CPU_step_gamma10000}

\begin{scriptsize}
\begin{center}
\renewcommand{\arraystretch}{1.3}
{\begin{tabular}{|c||c||c|c|c|c|c|c|c|}
\hline
\multicolumn{1}{|c||}{} & \multicolumn{1}{|c||}{} & \multicolumn{7}{c|}{$\beta$} \\
\cline{3-9}
$\nu$ & $l$ & $10^{-1}$ & $10^{-2}$ & $10^{-3}$ & $10^{-4}$ & $10^{-5}$ & $10^{-6}$ & $10^{-8}$ \\
\hline\hline
\multirow{3}{*}{$\frac{1}{100}$} & 
$1$ & 11.1 & 3.6 & 3.3 & 3.4 & 3.8 & 3.5 & 4.2 \\
\cline{3-9}
& $2$ & 54.1 & 46.4 & 58.4 & 70.5 & 106 & 141 & 282 \\
\cline{3-9}
& $3$ & 361 & 358 & 359 & 359 & 488 & 820 & 1616 \\
\hline
\multirow{3}{*}{$\frac{1}{250}$} & 
$1$ & 3.3 & 3.5 & 3.5 & 3.5 & 3.9 & 3.6 & 4.3 \\
\cline{3-9}
& $2$ & 46.8 & 46.9 & 70.4 & 70.9 & 118 & 141 & 272 \\
\cline{3-9}
& $3$ & 295 & 360 & 362 & 430 & 629 & 1090 & 1219 \\
\hline
\multirow{3}{*}{$\frac{1}{500}$} & 
$1$ & 3.3 & 3.6 & 3.4 & 3.6 & 3.8 & 3.6 & 4.2 \\
\cline{3-9}
& $2$ & 59.7 & 59 & 59 & 94.1 & 106 & 106 & 284 \\
\cline{3-9}
& $3$ & 362 & 363 & 362 & 493 & 625 & 825 & 1416 \\
\hline
\end{tabular}}
\end{center}
\end{scriptsize}
\end{table}

\section{Conclusions}\label{sec_5}
In this work, an augmented Lagrangian preconditioner has been proposed
	for use within an inexact Newton iteration applied to the control
	of the stationary Navier--Stokes equations. We compared the proposed
	preconditioner with a block pressure convection--diffusion preconditioner,
	which was derived for the Picard linearization of the
	Navier--Stokes control problem. Numerical results showed the good robustness
	of the augmented Lagrangian approach with respect to the mesh size, the
	regularization parameter, and the viscosity of the fluid, when employing
	an exact solver for the momentum blocks. The augmented Lagrangian
	preconditioner has been shown to be more robust than the (exact) block pressure
	convection--diffusion preconditioner when solving an inexact Newton
	linearization of the control problem considered here. Most importantly, the
	robustness of the approach is preserved even if the momentum block is solved
	inexactly, through the action of a suitable multigrid approach. Future research will focus
	on an eigenvalue analysis of the proposed strategy, its extension
	to the control of time-dependent Navier--Stokes equations,
	and the application of the augmented Lagrangian preconditioner to more
	complicated PDE-constrained optimization problems.

\section*{Acknowledgements}
The authors gratefully acknowledge two anonymous referees for their valuable comments. SL and PEF acknowledge Pablo Brubeck Martinez for useful discussions. PEF was funded by the Engineering and Physical Sciences Research Council [grant numbers EP/R029423/1 and EP/W026163/1]. SL and MB are members of Gruppo Nazionale di Calcolo Scientifico (GNCS) of the Istituto Nazionale di Alta Matematica (INdAM).

%
%
%

\end{document}